\documentclass[12pt]{amsproc}

\usepackage{amsmath}
\usepackage{amsthm}
\usepackage{amscd}
\usepackage{amssymb}
\usepackage{url}
\usepackage{nicefrac}
\usepackage[pdftex]{graphicx}%
\usepackage[outercaption]{sidecap}
\usepackage{amsfonts,epigraph}%
\usepackage{booktabs} 
\usepackage[all,cmtip]{xy}
\newtheorem{theorem}{Theorem}[section]
\newtheorem{corollary}[theorem]{Corollary}
\newtheorem{lemma}[theorem]{Lemma}
\newtheorem{proposition}[theorem]{Proposition}

\theoremstyle{definition}
\newtheorem{definition}[theorem]{Definition}
\newtheorem{example}[theorem]{Example}

\newtheorem{choices-notations}[theorem]{Choices and Notations}
\newtheorem{notation}[theorem]{Notation}

\newtheorem{rem}[theorem]{Remark}

\theoremstyle{remark}

\newcommand{\exendproof}{\renewcommand{\qed}{\relax}\end{proof}}

\newsavebox{\SmallMathBox}

\DeclareRobustCommand*{\nicefrac}[2]{\ifmmode\mathnicefrac{#1}
{ #2}%
  \else\textnicefrac{#1}{#2}\fi}
\newcommand*{\textnicefrac}[2]{\check@mathfonts%
\mbox{\raisebox{.5ex}{\fontsize\sf@size\z@\selectfont#1}\kern-.
1em%
/\kern-.1em\raisebox{- .25ex}{\fontsize\sf@size\z@\selectfont#2} }}
\newcommand*{\mathnicefrac}[2]{%
  \mathchoice
    {\m@fr@c{\scriptstyle}{#1}{#2}}
    {\m@fr@c{\scriptstyle}{#1}{#2}}
    {\m@fr@c{\scriptscriptstyle}{#1}{#2}}
    {\m@fr@c{\scriptscriptstyle}{#1}{#2}}}

\def\lla{\langle}

\def\noi{\noindent}

\def\rra{\rangle}
\def\sqm1{\sqrt{-1}}

\def\tand{\mbox{\ \rm  and }}

\def\too{\longrightarrow}

\def\wh{\widehat}
\def\wt{\widetilde}

\def\={\cong}
\def\>{\supset}
\def\<{\subset}

\def\12{\frac{1}{2}}
\def\0{^{\circ}}

\newcommand{\N}{\mathbb N}
\newcommand{\Z}{\mathbb Z}
\newcommand{\R}{\mathbb R}

\newcommand{\C}{\mathbb C}
\newcommand{\OO}{\mathrm O}

\newcommand{\U}{\mathrm U}
\newcommand{\mi}{\mathrm i}

\def\CC{{\mathbb C}}

\def\KK{{\mathbb K}}

\def\NN{{\mathbb N}}

\def\RR{{\mathbb R}}

\def\UU{{\mathbb U}}
\def\ZZ{{\mathbb Z}}

\def\Bb{{\mathcal B}}

\def\Ff{{\mathcal F}}

\def\Ll{{\mathcal L}}

\def\Rr{{\mathcal R}}
\def\Ss{{\mathcal S}}

\def\C{\CC}

\def\la{\lambda}

\def\N{\NN}

\def\R{\RR}

\def\w{\omega}

\def\Z{\ZZ}

 \DeclareMathOperator{\dist}{dist}

\DeclareMathOperator{\dom}{dom}

  \DeclareMathOperator{\Graph}{graph}
 
\DeclareMathOperator{\Hom}{Hom}

 \DeclareMathOperator{\image}{im}
\DeclareMathOperator{\Index}{index}

 \DeclareMathOperator{\ran}{im} 
  
 \DeclareMathOperator{\romS}{S} 

 \DeclareMathOperator{\Sp}{Sp} \DeclareMathOperator{\Span}{span}

\begin{document}

\title[Skew-adjoint linear relations]{Skew-adjoint linear relations between Banach spaces}
%
%
%
%


\author{Hanchen Li}
\address{Information Technology Operations Center, Bank of China,
	Beijing 10094, P. R. China} \email{hhlhc@126.com
	}

\author{Chaofeng Zhu}
\address{Chern Institute of Mathematics and Laboratory of Pure Mathematics and Combinatorics (LPMC), Nankai University,
	Tianjin 300071, P. R. China, https://orcid.org/0000-0003-4600-4253} \email{zhucf@nankai.edu.cn}
\thanks{Corresponding author: C. Zhu [\texttt{zhucf@nankai.edu.cn}]\\	
	Chaofeng Zhu is supported by National Key R\&D Program of China (2020YFA0713300), NSFC Grants (12571066, 11971245), Nankai Zhide Foundation and Nankai University.}

\date{today}




\begin{abstract}
In this paper, we prove the stability theorems for the isotropic perturbations of maximal isotropic subspaces in symplectic Banach spaces. Then we prove a stability theorem for the mod $2$ dimensions of kernel of skew-adjoint linear Fredholm relations between real Banach spaces with index $0$. Finally we gives the exactly two path components of the set of skew-adjoint linear Fredholm relations between real Banach spaces with indices $0$.
\end{abstract}

\maketitle


\section{Introduction}\label{l:introduction-stability-symplectic}

M. F. Atiyah and I. M. Singer in \cite{AS69} considered the set of all skew-adjoint Fredholm operators ${\mathcal F}^1(X)$ on a real separable Hilbert space $X$ and showed that it is a classifying space for (Atiyah's) real K-theory. As a consequence, the homotopy groups of ${\mathcal F}^1(X)$ are $8$-periodic, with a fundamental group that is cyclic of order $2$.
A.L. Carey, J. Phillips, H. Schulz-Baldes in \cite{CaPhSch19} constructed a homotopy invariant for paths in ${\mathcal F}^1(X)$ with invertible endpoints which provides an explicit isomorphism between $\pi_1({\mathcal F}^1(X))$ and ${\mathbb Z}_2$. The authors call it ${\mathbb Z}_2$-valued spectral flow, as formally it shares properties with the (${\mathbb Z}$-valued) spectral flow for paths of selfadjoint Fredholm operators (see M. F. Atiyah, V. K. Patodi and I. M. Singer\cite{APS76}). Recently, N. Doll\cite{Doll24} defined the oriented flow for paths in ${\mathcal F}^1(X)$ with odd dimensional kernels.

In the Hilbert space setting, one of the most powerful tool of studying selfadjoint operators and skew-adjoint operators is the spectral theory. By the idea of symplectic geometry, B. Booss and the second author \cite[Definition 4.2.1]{BoZh18} was able to generalize the notion of selfadjoint operators and skew-adjoint operators to the linear relations between two real (complex) linear vector spaces with a non-degenerate sesquilinear form. In the Banach space setting, the classical way of studying selfadjoint operators and skew-adjoint operators fails. Our paper is an attempt to the study of selfadjoint linear relations and skew-adjoint linear relations between Banach spaces. Then it is possible to define the oriented flow in the Banach space situation. Note that our definition of skew symmetric operators and skew-adjoint operators in Banach spaces (see Definition \ref{d:selfadjoint} below) is different from the one given by N. Ghoussoub \cite[Definition 4.1]{Gh09}.

Since the selfadjoint operators form the most important class of operators that appear in applications, the stability of selfadjointness is one of our important problems. In \cite[Section V.4]{Ka95}, T. Kato had given a nice introduction to this topic. 
Let $X$ be a Hilbert space. A linear operator $A\colon X\supset\dom(A)\to X$ is symmetric (or selfadjoint) if and only if the graph $\Graph(A)$ of $A$ is an isotropic (or Lagrangian) linear subspace of $X\times X$.

In this paper, we shall generalize the theory to the {\em symplectic} setting. Note that the proofs in \cite[Section V.4]{Ka95} do not apply to the symplectic case. 

Throughout this paper, we denote by $\KK$ the field of real numbers or complex numbers. We denote by $\N$, $\Z$, $\R$ and $\C$ the sets of all natural, integral, real and complex numbers respectively. By $\UU$ we denote the unit circle in the complex plane. We denote by $\mi$ the imaginary unit. Denote by $I_X$ the identity map on a set $X$. If there is no confusion, we will omit the subscript $X$. For a Banach space $X$ over $\KK$, we denote by $\Bb(X)$ the set of bounded linear operators. As in \cite[Sectioin III.1.4]{Ka95}, the {\em adjoint space} $X^*$ of a Banach space is defined to be the set of bounded semi-linear forms on $X$. We equip the the set of closed linear subspaces $\Ss(X)$ with {\em gap distance} $\hat\delta$ (\cite[Section IV.2.1]{Ka95}). Let $X,Y$ be two Banach spaces. A bounded sesquilinear form $\Omega\colon X\times Y\to \KK$ define two bounded linear maps $L_\Omega\colon X\to Y^*$ and $R_\Omega\colon Y\to X^*$ by 
\begin{align}\label{e:form-operator}
	\Omega(x,y)=(L_\Omega x)y=\overline{(R_\Omega y)x}\quad\text{ for all }x\in X\text{ and }y\in Y.
\end{align}

Let $(X,\omega)$ be a symplectic vector space over $\KK$ with a maximal isotropic subspace. Then $h_{\lambda}$ denotes the sign of $\lambda$ defined by Definition \ref{d:sign-max-isotroic} below. By definition, a maximal isotropic subspace $\lambda$ is Lagrangian if and only if $h_\lambda=0$.

We have the following stability results of maximal isotropic subspaces in symplectic Banach spaces.

\begin{theorem}\label{t:stable-max-isotropic}
	Let $(X,\omega_0)$ be a symplectic Banach space with a maximal isotropic subspace $\lambda$. We assume that $\gamma_\lambda>0$ holds if $h_\lambda\ne 0$.
	Let $(X,\omega)$ be a symplectic Banach space with a closed isotropic subspace of $\mu$. 
	We set 
	\begin{align*}		
		l(\mu,\omega)\ :&=			    
		\left(|h_\lambda|\gamma_\lambda^{-1}\left(\|L_\omega-L_{\omega_0}\|+\|L_{\omega_0}\|\delta(\mu^\omega,\lambda^{\omega_0})(2+\delta(\mu^\omega,\lambda^{\omega_0}))\right)\right)^{\frac{1}{2}}. 				
	\end{align*}
	Here $h_\lambda$ and $\gamma_\lambda$ are defined by Definition \ref{d:sign-max-isotroic} below, and $l(\mu,\omega)=0$ if $h_\lambda=0$.
	
	Assume that there hold that 
	\begin{align}\label{e:stable-maximal-isotropic}
		\delta(\mu^\omega,\lambda^{\omega_0})+l(\mu,\omega)
		<\frac{1-\delta(\lambda,\mu)}{1+\delta(\lambda,\mu)}.
	\end{align}
	Then $\mu$ is a maximal isotropic subspace of $(X,\omega)$ with $h_\mu\in\{h_\lambda,0\}$.
\end{theorem}

We can use \cite[Lemma 1]{HeKa70} to prove Theorem \ref{t:stable-max-isotropic} if $h_\lambda=0$. To prove the general case, we need to study the augmented Morse index for the perturbation of a positive semi-definite symmetric form with possible different domain in a fixed Banach space. Note that in general one can not make the domains fixed.

By the study of the Morse index for the perturbation of a positive semi-definite symmetric form and Theorem \ref{t:stable-max-isotropic}, we have the following result.

\begin{theorem}\label{t:family-max-isotropic}
	Let $X$ be a Banach space with a continuous family of symplectic structure $\{\omega(s)\}_{s\in J}$, where $J$ is a connected topological space. Let $\{\lambda(s)\}_{s\in J}$ be a continuous family of closed subspace of $X$ such that $\lambda(s)\subset \lambda(s)^{\omega(s)}$ holds for each $s\in J$. Assume that $\{\lambda(s)^{\omega(s)}\}_{s\in J}$ and $\{\lambda(s)^{\omega(s)\omega(s)}\}_{s\in J}$ are also continuous families of closed subspace of $X$.
	Assume that there is an $s_0\in J$ such that $\lambda(s_0)$ is a maximal isotropic subspace of $(X,\omega(s_0))$ satisfying $\gamma_{\lambda(s_0)}>0$ if $h_{\lambda(s_0)}\ne 0$, and each $s\in J$ such that $\lambda(s)$ is a maximal isotropic subspace of $(X,\omega(s))$ satisfying $\gamma_{\lambda(s)}>0$ if $h_{\lambda(s)}\ne 0$. Then $\lambda(s)$ is a maximal isotropic subspace of $(X,\omega(s))$ and $h_{\lambda(s)}=h_{\lambda(s_0)}$ holds for each $s\in J$. 
\end{theorem}

In strong symplectic Banach spaces, some conditions in stability theorems of maximal isotropic subspaces hold. Then these stability theorems can be simplified(see Section \ref{ss:stable-max-isotropic-strong} below).

In 2023, J.M. F. Castillo, W. Cuellar
M. Gonz\'{a}lez and R. Pino \cite[Proposition 8.11]{CCGP23} showed that a strong real symplectic Hilbert space $(X,\omega)$ can not be a codimension $1$ linear subspace of a strong real symplectic Hilbert space $(Y,\tilde\omega)$ with compact $L_{\tilde\omega|_X}-L_\omega$. After their proof, they wrote:

{\em "We do not know whether Proposition 8.11 holds for general Banach spaces."}

Note that \cite[Proposition 8.11]{CCGP23} is a special case of the stability theorem for the mod 2 dimensions of kernel of bounded skew-adjoint operators of M. F. Atiyah and I. M. Singer \cite[p.6]{AS69} in real Hilbert spaces. In this paper, we shall solve the problem and generalize it to a stability theorem for the mod 2 dimensions of kernel and indeterminacy in real Banach spaces. Here we require that $L_\Omega$ is an injective bounded operator only.

\begin{theorem}\label{t:stable-mod-2-index}
	Let $X$ and $Y$ be two real Banach space, and $\Omega_0\colon X\times Y\to\R$ be a bounded non-degenerate bilinear form. 
	
	Let $T:X\leadsto Y$ be a skew-adjoint linear Fredholm relation with respect to $\Omega_0$ with index $0$.
	Then there is a $\delta>0$ such that for each skew symmetric closed linear relation $S:X\leadsto Y$ with respect to a bounded non-degenerate bilinear form $\Omega$ with   $\hat\delta(S,T)+\|L_\Omega-L_{\Omega_0}\|<\delta$, there hold that $S$ is skew-adjoint linear Fredholm relation with index $0$, and  
	\begin{align}\label{e:stable-mod-2-ker}
		\dim\ker S=\dim\ker T\mod 2.
	\end{align}
\end{theorem}

\begin{theorem}\label{t:family-mod-2-index}
	Let $X$ and $Y$ be two real Banach space with a continuous family of bounded non-degenerate bilinear forms  $\Omega(s)\colon X\times Y\to\R$, i.e. $\{L_{\Omega(s)}\in\Bb(X,Y^*)\}_{s\in J}$ is a continuous family of bounded injective operators, where $J$ is a connected topological space. 
	
	Let $\{T(s):X\leadsto Y\}_{s\in J}$ be a continuous family of linear Fredholm relations with index $0$ such that $T(s)$ is skew-adjoint with respect to $\Omega(s)$ for each $s\in J$.
	Then $\dim\ker T(s)\mod 2$ is a constant in $\ZZ_2$ for $s\in J$.
\end{theorem}

The following corollary answers the question raised in \cite{CCGP23} affirmatively.

\begin{corollary}\label{c:real-symplectic-extension-even}
	Let $(X,\omega)$ and $(X\oplus\R^n,\tilde\omega)$ be two real  symplectic Banach space, where $n$ is a non-negative integer. 
	Assume that $L_{\tilde\omega|_X}-L_\omega$ is compact. Then $(X,\omega)$ is a strong symplectic Banach space if and only if $(X\oplus\R^n,\tilde\omega)$ is. In the affirmative case, $n$ is even.
\end{corollary}

The path components of the set of skew-adjoint Fredholm relations between real Banach spaces with index $0$ is given as follows.

\begin{theorem}\label{t:path-components-real-skew}
	Let $X$ and $Y$ be two real Banach space, and $\Omega\colon X\times Y\to\R$ be a bounded non-degenerate bilinear form. 
	
	Let $\Ff\Rr_0^{-1,sa}(\Omega)$ be the set of real skew-adjoint linear Fredholm relations with index $0$ with respect to $\Omega$. For $a\in\{0,1\}$, we define 
	\begin{align}\label{e:skew-adjoint-a}
		\Ff\Rr_{a;0}^{-1,sa}(\Omega)\ :&=\{T\in \Ff\Rr_0^{-1,sa}(\Omega);\;\dim\ker T=a\mod 2\}.
	\end{align}
    Then the space $\Ff\Rr_0^{-1,sa}(\Omega)$ has exactly two path components $\Ff\Rr_{0;0}^{-1,sa}(\Omega)$ and $\Ff\Rr_{1;0}^{-1,sa}(\Omega)$.
\end{theorem}

In the finite dimensional case, we can always give a Hilbert structure on $X$. By Proposition \ref{p:real-skew-adjoint-Hilbert} below, there is a homoemorphism $\varphi_\Omega\colon\mathrm{O}(X)\to\Ff\Rr_0^{-1,sa}(\Omega)=\Rr^{-1,sa}(\Omega)$. Moreover, we have 
\begin{align*}
	\ker(T)=\ker(\varphi_\Omega^{-1}(T)-I_X)
\end{align*}
for each $T\in\Rr^{-1,sa}(\Omega)$. Thus Theorem \ref{t:stable-mod-2-index} and Theorem \ref{t:path-components-real-skew} hold in this case. 

In the general Banach space situation, we use the symplectic reduction method and reduce the problems to the finite dimensional case.

The paper is organized as follows. In Section 1, we explain the motivation of the research and state the main results of the paper. In Section 2, we prove Theorem \ref{t:stable-max-isotropic} and Theorem \ref{t:family-max-isotropic}. In Section 3, we study the stability theorems with relative bounded perturbations. In Section 4, we prove Theorem \ref{t:stable-mod-2-index} and Theorem \ref{t:path-components-real-skew}. 

We would like to thank the referees of this paper for their critical reading and
very helpful comments and suggestions.

\section{Stability theorems in symplectic Banach space}\label{s:stable-symplectic}

\subsection{Fredholm pairs of closed linear subspaces}\label{ss:Fredholm}

Let $X$ be a Banach space. We denote by $\Ss(X)$ the set of closed linear subspaces of $X$ and $\Ss^c(X)$ the set of complemented closed linear subspaces of $X$ respectively. Let $M,N$ be two closed linear subspaces of $X$, i.e., $M,N\in\Ss(X)$.
Denote by $\romS_M$ the unit
sphere of $M$. We recall three common definitions of distances in $\Ss(X)$ (see also
\cite[Sections IV.2.1 and IV.4.1]{Ka95}):
\begin{itemize}
	
	\item  the \textit{Hausdorff metric} $\hat d$;
	
	\item the \textit{aperture (gap distance)} $\hat\delta$, that is not a metric since it does not in general
	satisfy the triangle inequality, but defines the same topology as the metric $\hat d$, called
	\textit{gap topology}, and is easier to estimate than $\hat d$; and
	
	\item the \textit{angular distance (minimum gap)} $\wh\gamma$, that is useful in our estimates, though not
	defining any suitable topology.
	
\end{itemize}

\begin{definition}[The gap between closed linear subspaces]\label{d:closed-distance}	
	(a) We set \begin{multline*}
		\hat d(M,N)\ =\
		d(\romS_M,\romS_N)\\
		:= \begin{cases} \max\left\{\begin{matrix}\sup\limits_{u\in \romS_M}\dist(u,\romS_N),\\ \sup\limits_{u\in
					\romS_N}\dist(u,\romS_M)\end{matrix}\right\},&\text{ if both $M\ne 0$ and $N\ne 0$},\\
			0,&\text{ if $M=N=0$},\\
			2,&\text{ if either $M= 0$ and $N\ne 0$ or vice versa}.
		\end{cases}
	\end{multline*}
	\newline (b) We set
	\begin{eqnarray*}
		\delta(M,N)\ &:=&\ \begin{cases}\sup\limits_{u\in \romS_M}\dist(u,N),& \text{if $M\ne\{0\}$},\\
			0,& \text{if $M=\{0\}$},\end{cases}\\
		\hat\delta(M,N)\ &:=&\ \max\{\delta(M,N),\delta(N,M)\}.
	\end{eqnarray*}
	$\hat\delta(M,N)$ is called the {\em gap} between $M$ and $N$.
	\newline (c) We set
	\begin{align*}
		\gamma(M,N)\ :=&\ \begin{cases} \inf\limits_{u\in M\setminus N}\frac{\dist(u,N)}{\dist(u,M\cap N)}\ (\le
			1), & \text{if $M\not\subset N$},\\
			1, & \text{if $M\subset N$},\end{cases}\\
		\hat\gamma(M,N)\ :=&\ \min\{\gamma(M,N),\gamma(N,M)\}.
	\end{align*}
	$\hat\gamma(M,N)$ is called the {\em minimum gap} between $M$ and $N$. If $M\cap N=\{0\}$, we have
	\[\gamma(M,N)\ =\ \inf\limits_{u\in \romS_M}\dist(u,N).\]
\end{definition}

\begin{rem}\label{r:delta-between-linear-subspaces}
	Let $X$ be a Banach space with two linear subspaces $M, N$. Then we have $\hat d(M,N)=\hat d(\bar M,\bar N)$, $\delta(M,N)=\delta(\bar M,\bar N)$ and $\hat\delta(M,N)=\hat\delta(\bar M,\bar N)$.
	However, $\hat d$ is not a metric on the space of all linear subspaces of $X$ if $\dim X=+\infty$.
\end{rem}

In this paper we shall impose the gap topology on the space $\Ss(X)$ of all closed linear
subspaces of a Banach space $X$ and its subset $\Ss^c(X)$ of complemented subspaces.

\smallskip

\begin{definition}\label{d:fredholm-pair}
	(a) The space of (algebraic) \emph{Fredholm pairs} of linear subspaces of a vector space $X$ is
	defined by
	\begin{equation}\label{e:fp-alg}
		\Ff^{2}_{\operatorname{alg}}(X)\ : =\{(M,N);\;  \dim (M\cap N)  <+\infty,\; \dim
		X/(M+N)<+\infty\}
	\end{equation}
	with
	\begin{equation}\label{e:fp-index}
		\Index(M,N;X)=\Index(M,N)\ :=\dim(M\cap N) - \dim X/(M+N).
	\end{equation}
	
	\noi (b) In a Banach space $X$, the space of (topological) \emph{Fredholm pairs} is defined by
	\begin{multline}\label{e:fp}
		\Ff^{2}(X)\ :=\{(M,N)\in\Ff^2_{\operatorname{alg}}(X);\; M,N, \tand M+N \subset X \text{ closed}\}.
	\end{multline}
	A pair $(M,N)$ of closed subspaces is called {\em semi-Fredholm} if $M+N$ is closed, and at least
	one of $\dim (M\cap N)$ and $\dim\left(X/(M+N)\right)$ is finite.
	
	\noi (c) Let $X$ be a Banach space, $M\in\Ss(X)$ and $k\in\Z$. We define
	\begin{align}
		\label{e:fp-M}\Ff_M(X)\ :&=\{N\in\Ss(X);\;(M,N)\in\Ff^2(X)\},\\
		\label{e:fp-kM}\Ff_{k,M}(X)\ :&=\{N\in\Ss(X);\;(M,N)\in\Ff^2(X),\Index(M,N)=k\}.
	\end{align}
\end{definition}

\subsection{Symplectic Banach spaces}\label{ss:symplectic}
Firstly we recall the basic concepts and properties of symplectic functional analysis.

Let $X$ be a vector space over $\KK$, and $\Omega\colon X\times X\to\KK$ be a symmetric sesquilinear form or a skew symmetric sesquilinear form.
The {\em annihilator} of a subspace $\lambda$ of $X$ is defined by
\begin{equation}\label{e:annihilator}
	\lambda^\Omega\ :=\{y\in X;\;\Omega(x,y)=0, \; \forall x\in\lambda \}.
\end{equation}
The form $\Omega$ is called {\em non-degenerate}, if $X^{\Omega}=\{0\}$.

\begin{definition}\label{d:symplectic-space}
	Let $X$ be a vector space over $\KK$.%
	\newline (a) A mapping
	\[
	\omega\colon X\times X\too \KK
	\]
	is called a {\em symplectic form} on $X$, if it is a
	non-degenerate skew symmetric sesquilinear form.
	
	Then we call $(X,\w)$ a {\em symplectic vector space}.
	\newline (b) Let $X$ be a Banach space and $(X,\omega)$ a symplectic vector space. $(X,\omega)$
	is called {\em (weak) symplectic Banach space} if $\omega$ is bounded, i.e., $|\w(x,y)| \leq
	C\|x\| \|y\|$ for all $x,y\in X$. $(X,\omega)$
	is called {\em strong symplectic Banach space} if $\omega$ is bounded and $L_\omega^{-1}\in\Bb(X^*,X)$.
	\newline (c) Let $(X_j,\omega_j)$ be two symplectic vector spaces. A linear map $L\colon X_1\to X_2$ is called {\em symplectic} if $L^*\omega_2=\omega_1$, i.e., $\omega_2(Lx,Ly)=\omega_1(x,y)$ for each $x,y\in X_1$. 
	The {\em symplectic group} $\Sp(X,\omega)$ consists of all surjective symplectic linear maps of a symplectic vector space $(X,\omega)$. We write $\Sp(X)\ :=\Sp(X,\omega)$ if there is no confusion.
	\newline (c) A subspace ${\la}$ is called {\em symplectic},
	{\em isotropic}, {\em co-isotropic}, or {\em Lagrangian} if
	\[
	\la\cap{\la}^{\w}\ =\ \{0\}\,,\quad{\la} \,\subset\, {\la}^{\w}\,,\quad {\la}\,\supset\, {\la}^{\w}\,,\quad
	{\la}\,\ =\ \, {\la}^{\w}\,,
	\]
	respectively. An isotropic subspace ${\la}$ is called {\em maximal isotropic}, if each isotropic subspace $\mu\supset\lambda$ satisfies $\mu=\lambda$.
	\newline (d) The {\em Lagrangian Grassmannian}
	$\Ll(X,\w)$ consists of all Lagrangian subspaces of
	$(X,\w)$. We write $\Ll(X)\ :=\Ll(X,\w)$ if there is no confusion.
	\newline (d) The {\em complemented Lagrangian Grassmannian}
	$\Ll^c(X,\w)$ is defined by 
	\begin{equation}\label{e:complemented-Lagrangian}
		\Ll^c(X,\w)\ :=\{\lambda\in\Ll(X,\w);\;X=\lambda\oplus\mu\text{ for some }\mu\in\Ll(X,\w)\}.
	\end{equation}
	We write $\Ll^c(X)\ :=\Ll^c(X,\w)$ if there is no confusion.
\end{definition}

For a symplectic Banach space $(X,\omega)$, the induced operator $L_\omega\colon X\to X^*$ defined by \eqref{e:form-operator} is skew-adjoint in the sense that $L_\omega^*|_X=-L_\omega$.

A {\em linear relation} $A\colon X\rightsquigarrow Y$
between two linear space $X$ and $Y$ is a linear subspace of $X\times Y$. For two linear relations $A\colon X\rightsquigarrow Y$ and $B\colon Y\rightsquigarrow Z$, the {\em composite} $B\circ A$ is defined by
\begin{align}\label{e:composite}
	B\circ A=\{(x,z)\in X\times Z;\ (x,y)\in A\text{ and }(y,z)\in B\text{ for some }y\in Y\}.
\end{align}
A linear operator can be viewed as a linear relation by taking its graph. For details, we refer to R. Cross \cite{Cr98} and G. W. Whitehead \cite[B.2]{GWh78}.

\begin{notation}\label{n:confusing-notions} 
	Let $X$ and $Y$ be two linear spaces over $\KK$, $A$ and $B$ be two linear subspaces of $X\times Y$, and $a\in\KK$ be a number. We shall use different notations to distinguish the following notions if there is confusion:
	\begin{itemize}
		\item as linear subspaces, we have
		\begin{align*}
			A+B=&\{u+v;\;u\in A,v\in B\},\\
			aA=&\{au;\;u\in A\};
		\end{align*}
		\item as linear relations, we have
		\begin{align*}
			A\hat+B=&\{(x,y+z);\;(x,y)\in A,(x,z)\in B\},\\
			a\circ A=&\{(x,ay);\;(x,y)\in A\}.
		\end{align*}
	\end{itemize}	
\end{notation}

Let $X$ be a Banach space and $M\in\Ss(X)$ be a closed subspace. We define 
\begin{align}\label{e:M-bot}
	M^\bot\ :=\{f\in X^*;\; f(x)=0,\text{ for all } x\in M\}.
\end{align}

Let $X$ and $Y$ be two vector spaces over $\KK$. Given a non-degenerate sesquilinear form $\Omega\colon X\times Y\to\KK$, there is a natural symplectic structure $\omega$ on $Z\ :=X\times Y$ defined by (cf. \cite{Ben72}, \cite[Definition 4.2.1]{BoZh18}, \cite{Ek90})
\begin{align}\label{e:natural-symplectic-structure}
	\omega((x_1,y_1),(x_2,y_2))\ :=\Omega(x_1,y_2)-\overline{\Omega(x_2,y_1)}
\end{align}
for all $x_1,x_2\in X$ and $y_1,y_2\in Y$. For each subsets $\lambda$ of $X$ and $\mu$ 
of $Y$, the {\em annihilator} 
$\lambda^{\Omega,r}$ and $\mu^{\Omega,l}$ (cf. \cite[Definition 1.1.1.a]{BoZh18}) are defined by 
\begin{align}\label{e:right-side-annihilator}
	\lambda^{\Omega,r}\ :=&\{y\in Y;\ \Omega(x,y)=0\text{ for all }x\in\lambda\},\\
	\label{e:left-side-annihilator}
	\mu^{\Omega,l}\ :=&\{x\in X;\ \Omega(x,y)=0\text{ for all }y\in\mu\}.
\end{align}
Then there are two nature Lagrangian subspaces $X\times\{0\}$ and $\{0\}\times Y$ of $(Z,\omega)$. If $X, Y$ are Banach spaces and $\Omega$ is bounded, by \eqref{e:form-operator}, the induced operator of $\omega$ is given by
\begin{align}\label{e:L-omega}
	L_\omega=\left(\begin{array}{cc}
		0&-R_\Omega  \\
		L_\Omega&0 
	\end{array}\right).
\end{align}

Let $X$ be locally convex space. We define $\Omega_X\colon X\times X^*\to\KK$ by $\Omega_X(x,f)=\overline{f(x)}$ for all $x\in X$ and $f\in X^*$. By Hahn-Banach theorem, $\Omega_X$ is nondegenerate and continuous. Then $(X\times X^*,\omega)$ is a symplectic locally convex space. 

\begin{definition}(cf. \cite[Definition 4.2.1]{BoZh18})\label{d:selfadjoint}
	Let $X$ and $Y$ be two vector spaces over $\KK$. Let $\Omega\colon X\times Y\to\KK$ be a nondegenerate sesquilinear form and $h\in\KK$ be a number. Let $\omega$ be the associated symplectic form of $\Omega$ defined by \eqref{e:natural-symplectic-structure}.
	\begin{itemize}
		\item [(a)] Let $A\colon X\leadsto Y$ be a linear relation. The {\em $\Omega$ adjoint} of $A$ is defined to be $A^{\Omega}\ :=A^{\omega}$. We call $A$ {\em $h$ symmetric with respect to $\Omega$} if $h\circ A\subset A^{\Omega}$, and {\em $h$ selfadjoint with respect to $\Omega$} if $h\circ A= A^{\Omega}$, respectively. If $h=1$, we call $A$ {\em symmetric with respect to $\Omega$} and {\em selfadjoint with respect to $\Omega$} respectively. If $h=-1$, we call $A$ {\em skew symmetric with respect to $\Omega$} and {skew-adjoint with respect to $\Omega$} respectively. 
		We call $A$ a {\em maximal $h$ symmetric linear relation with respect to $\Omega$}, if $A$ is an $h$ symmetric linear relation with respect to $\Omega$ and maximal in the set of all $h$ symmetric linear relations with respect to $\Omega$. 
		\item [(b)] We denote by $\Rr^{h,sa}(\Omega)$ the set of $h$ selfadjoint relations between $X$ and $Y$ with respect to $\Omega$. We denote by $\Rr^{h,sa}(X)\ : =\Rr^{h,sa}(\Omega_X)$. For $k\in\Z$, we define
		\begin{align*}
			\Ff\Rr^{h,sa}(\Omega)\ :&=\{A\in \Rr^{h,sa}(\Omega);\; A\text{ is Fredholm}\},\\
			\Ff\Rr_k^{h,sa}(\Omega)\ :&=\{A\in \Ff\Rr^{h,sa}(\Omega);\; \Index A=k\},\\
			\Ff\Rr^{h,sa}(X)\ :&=\Ff\Rr^{h,sa}(\Omega_X),\\
			\Ff\Rr_k^{h,sa}(X)\ :&=\Ff\Rr_k^{h,sa}(\Omega_X).
		\end{align*}
		\item[(c)] Let $A\colon X\leadsto Y$ be a $h$-symmetric linear relation with $h\ne 0$. Define $Q_A\colon\dom (A)\times\dom(A)\to\KK$ by $Q_A(x,y)\ :=\Omega(x,Ay)$ for all $x,y\in\dom(A)$. Then $Q_A$ is a well-defined sesquilinear form. We call $Q_A$ the {\em associated form} of $A$.
	\end{itemize}	
\end{definition}

By definition, we have 
\begin{align}\label{e:A-Omega-trivial-subset}
	A^\Omega\supset(\ran A)^{\Omega,l}\times(\dom(A))^{\Omega,r}.
\end{align}

\begin{example}\label{ex:1-dimension-case}
	(a) We have    $\Rr^{-1,sa}(\R)=\{\R\times\{0\},\{0\}\times\R\}$.

	\noi(b) We have
	 \begin{align*}
	 	\Rr^{1,sa}(\Omega)\cap\Rr^{-1,sa}(\Omega)\supset\{X\times\{0\},\{0\}\times Y\}.
	 \end{align*}
\end{example}

The following lemma characterize the linear relations which are both symmetric and skew symmetric.

\begin{lemma}\label{l:symmetric-cap-skew-symmetric}
	Let $X$ and $Y$ be two vector spaces over $\KK$. Let $\Omega\colon X\times Y\to\KK$ be a non-degenerate sesquilinear form and $h_1,h_2\in\KK$ be two different numbers with $|h_1|=|h_2|=1$. Let $A\colon X\leadsto Y$ be a linear relation. Then the following hold.
	\begin{itemize}
		\item[(a)] $A$ is $h_1$ symmetric and $h_2$ symmetric with respect to $\Omega$ if and only if $A^\Omega\supset\dom(A)\times\ran A$, if and only if $\dom(A)\subset(\ran A)^{\Omega,l}$. 
		\item[(b)] Assume that $A$ is $h_1$ symmetric and $h_2$ symmetric with respect to $\Omega$. Then $A$ is a maximal $h_1$ symmetric linear relation with respect to $\Omega$ if and only if $A=\dom(A)\times\ran A$, $\dom(A)=(\ran A)^{\Omega,l}$ and $\ran A=(\dom(A))^{\Omega,r}$, if and only if $A$ is $h$ selfadjoint with respect to $\Omega$ for $|h|=1$. 
	\end{itemize}
\end{lemma}

\begin{proof}
	(a) (i) If $A^\Omega\supset\dom(A)\times\ran A$, we have $(h_1\circ A)\cup(h_2\circ A)\subset\dom(A)\times\ran A\subset A^\Omega$. Then $A$ is $h_1$ symmetric and $h_2$ symmetric with respect to $\Omega$.
	
	(ii) Assume that $(h_1\circ A)\cup(h_2\circ A)\subset A^\Omega$ holds. Let $(x,y)$ be in $A$. Then we have $(x,h_1y)\in A^\Omega$ and $(x,h_2y)\in A^\Omega$. Then we obtain $(x,0)\in A^\Omega$ and $(0,y)\in A^\Omega$. For each $y_1\in\ran A$, we have $(0,y_1)\in A^\Omega$. Thus we have 
	\begin{align*}
		\Omega(x,y_1)=\overline{\Omega(0,y)}=0.
	\end{align*}
	Then we have $A^\Omega\supset\dom(A)\times\ran A$ and $\dom(A)\subset(\ran A)^{\Omega,l}$. 
	
	(iii) Assume that $\dom(A)\subset(\ran A)^{\Omega,l}$ holds. Since $(h_1\circ A)\cup(h_2\circ A)\subset\dom(A)\times\ran A$, we have $A^\Omega\supset\dom(A)\times\ran A$. By Step a.i, we have $(h_1\circ A)\cup(h_2\circ A)\subset A^\Omega$. Then $A$ is $h_1$ symmetric and $h_2$ symmetric with respect to $\Omega$.
	
	\noi (b) (i) By (a), we have $h_1\circ A\subset\dom(A)\times\ran A\subset A^\Omega$, and $\dom(A)\times\ran A$, $\dom(A) \times(\dom(A))^{\Omega,r}$, and $(\ran A)^{\Omega,l}\times\ran A$ are three $h_1$ symmetric linear relations with respect to $\Omega$ containing $A$. 
	
	(ii) If $A$ is a maximal $h_1$ symmetric linear relation with respect to $\Omega$, by Step b.i, we have $A=\dom(A)\times\ran A=\dom(A) \times(\dom(A))^{\Omega,r}=(\ran A)^{\Omega,l}\times\ran A$. Then we have $A=\dom(A)\times\ran A$, $\dom(A)=(\ran A)^{\Omega,l}$ and $\ran A=(\dom(A))^{\Omega,r}$. Then $A$ is $h$ selfadjoint with respect to $\Omega$ for $|h|=1$.
	
	(iii) By Step b.ii, we have $A=\dom(A)\times\ran A=\dom(A) \times(\dom(A))^{\Omega,r}=(\ran A)^{\Omega,l}\times\ran A$. Then $A$ is a maximal $h_1$ symmetric linear relation with respect to $\Omega$.
\end{proof}

Let $X, Y$ be two Banach spaces. 
Let $A\colon X\leadsto Y$ be a linear relation. By \cite[(III.5.9)]{Ka95} we have 
\begin{align}\label{e:addjoint-relation}
	(-1\circ A^*)^{-1}=A^\perp.
\end{align}
Given a bounded nondegenerate sesquilinear form $\Omega\colon X\times Y\to\KK$, we have
\begin{align}\label{e:Omega-addjoint-relation}
	A^\Omega=(L_\omega)^{-1}A^\perp=(L_\omega)^{-1}(-1\circ A^*)^{-1}=(R_\Omega)^{-1}\circ A^*\circ L_\Omega,
\end{align}
where $L_\omega$ is defined by \eqref{e:L-omega}, $L_\Omega$ and $R_\Omega$ are defined by (\ref{e:form-operator}), and $(L_\omega)^{-1}A^\perp=(L_\omega)^{-1}(-1\circ A^*)^{-1}$ means the operator $(L_\omega)^{-1}$ acts on the linear subspace $A^\perp=(-1\circ A^*)^{-1}$.

Let $X$ be a Banach space. We denote by $\iota_X\colon X\to X^{**}$ the natural embedding and $\mi\ :=\sqrt{-1}$. Then we have $R_\Omega=L_\Omega^*\circ\iota_Y$,  $L_{\Omega_X}=\iota_X$ and $R_{\Omega_X}=I_{X^*}$. By \eqref{e:Omega-addjoint-relation} we have $A^{\Omega_X}=A^*\circ\iota_X$. Then $A$ is $h$ symmetric if and only if $h\circ A\subset A^*\circ\iota_X$, and $A$ is $h$ selfadjoint if and only if $h\circ A=A^*\circ\iota_X$. 
Assume that $a\in \KK=\C$ and $|a|=1$. Then $A$ is $a^2$ symmetric if and only if $a\circ A$ is symmetric, $A$ is $a^2$ selfadjoint if and only if $a\circ A$ is selfadjoint. Thus we have
\begin{align}\label{e:a2-sa}
	a\circ\Rr^{a^2,sa}(\Omega)=\Rr^{1,sa}(\Omega).
\end{align}

The number $|h|$ can only take values of $1$ in the non-trivial case. More precisely, we have the following result. 

\begin{lemma}\label{l:value-h}
	Let $X$ and $Y$ be two vector spaces over $\KK$. Let $\Omega\colon X\times Y\to\KK$ be a nondegenerate sesquilinear form and $h\in\KK$ be a number.
	Let $A\colon X\leadsto Y$ be an $h$ selfadjoint linear relation with respect to $\Omega$. Then we have the following.
	\begin{itemize}
		\item [(a)] There holds $A^{\Omega}=0\circ A$ if and only if $A=X\times \{0\}$.
		\item [(b)] There holds $|h|= 1$ if $A\ne X\times\{0\}$ and $A\ne \{0\}\times Y$.
	\end{itemize}
\end{lemma}

\begin{proof} 
	(a) If $A=X\times\{0\}$, we have $A^\Omega=A=0A$.
	 
	Assume that $A^{\Omega}=0\circ A$. By definition, we have 
	\begin{align*}
		\Omega(x_1,y_2)=\overline{\Omega(x_2,0)}=0
	\end{align*}
	for each $x_1\in\dom(A)$ and $(x_2,y_2)\in A$. Then we have $\dom(A)\subset(\ran A)^{\Omega,l}$. By \eqref{e:A-Omega-trivial-subset}, we have $A^\Omega\supset(\ran A)^{\Omega,l}\times(\dom(A))^{\Omega,r}$. Since $A^{\Omega}=0 A$ holds, we have $\dom(A)=(\ran A)^{\Omega,l}$ and $(\dom(A))^{\Omega,r}=\{0\}$. Then we have 
	\begin{align*}
		\ran A\subset((\ran A)^{\Omega,l})^{\Omega,r}=(\dom(A))^{\Omega,r}=\{0\}.
	\end{align*}
	Thus we have $\dom(A)=(\ran A)^{\Omega,l}=X$ and $A=X\times \{0\}$.
	\newline (b) (i) By (a), there holds $h\ne 0$.
	
	(ii) {\em There holds $|h|= 1$ if $A\not\subset X\times\{0\}$ and $A\not\subset \{0\}\times Y$.}
	
	Since $A^{\Omega}=h\circ A$ and $h\ne 0$, we have 
	\begin{align*}
		A\subset A^{\Omega\Omega}= (h A)^{\Omega}=\bar h\circ A^\Omega=|h|^2\circ A.
	\end{align*}
    Since $A\not\subset X\times\{0\}$ and $A\not\subset \{0\}\times Y$, there is an $x\in X\setminus\{0\}$ and a $y\in Y\setminus\{0\}$ such that $(x,y)\in A$. Then $0\notin Ax$ and $Ax$ is not a linear subspace of $Y$. Since $Ax=|h|^2Ax$, we have $|h|^2=1$ and $|h|=1$.
    
    (iii) Assume that $A\ne X\times\{0\}$ and $A\ne \{0\}\times Y$. By (a), we have $A^\Omega\ne 0\circ A$ and $(A^{-1})^\Omega\ne 0\circ A^{-1}$. Since $A^\Omega=h\circ A$ and $h\ne 0$, we have $A\not\subset X\times\{0\}$ and $A\not\subset \{0\}\times Y$. By Step b.ii, we have $|h|=1$.
\end{proof}

An $h$ symmetric linear relation ($|h|=1$) has the following property.

\begin{lemma}\label{l:ker-dom-symmetric}
	Let $X$ and $Y$ be two vector space over $\KK$. Let $\Omega\colon X\times Y\to\KK$ be a nondegenerate sesquilinear form with respect to $\Omega$. Let $h\in\KK$ be such that $|h|=1$. 
	Let $A\colon X\leadsto Y$ be an $h$ symmetric linear relation. Then we have 
	$\Omega(x,y)=0$ for $(x,y)\in(\ker A\times\image A)\cup(\dom(A)\times A0)$.
\end{lemma}

\begin{proof}
	Since $(x,y)\in A$ and $A\subset A^\Omega$, we have $(x,y)\in A^\Omega$. Since $x=0$ or $y=0$ holds, by \eqref{e:natural-symplectic-structure} we have $\Omega(x,y)=0$.	
\end{proof}

For an $h$ maximal symmetric linear relation ($|h|=1$), we can say more.

\begin{proposition}\label{p:ker-dom-max-symmetric}
	Let $X$ and $Y$ be two vector space over $\KK$. Let $\Omega\colon X\times Y\to\KK$ be a nondegenerate sesquilinear form. Let $h\in\KK$ be such that $|h|=1$. 
	Let $A\colon X\leadsto Y$ be a maximal $h$ symmetric linear relation with respect to $\Omega$. Then we have 
	\begin{align}\label{e:ker-dom}
			&\ker A=\ker A^\Omega=(\image A)^{\Omega,l}=\ker Q_A,\\
			\label{e:ker0-dom}
			&A0=A^\Omega0=(\dom(A))^{\Omega,r},\\
			\label{e:im-dom-1}
			&\image A^\Omega\subset(\ker A)^{\Omega,r}=((\image A)^{\Omega,l})^{\Omega,r},\\
			\label{e:im-dom-2} 
			&\dom A^\Omega\subset(A0)^{\Omega,l}=((\dom A)^{\Omega,r})^{\Omega,l}.
	\end{align}		
\end{proposition}

\begin{proof} Let $\omega$ be the symplectic structure defined by \eqref{e:natural-symplectic-structure}.%
	\newline 1. Let $x$ be in $(\image A)^{\Omega,l}$. By \eqref{e:natural-symplectic-structure}, we have
	$(x,0)\in A^\Omega=A^\omega$. Note that we have $\omega((x,0),(x,0))=0$. Set $V\ :=\Span\{(x,0)\}$ and $W\ :=A+V$. Note that $h\circ V=V$. Then we have 
	\begin{align*}
		W^\Omega=A^\Omega\cap V^\Omega
		\supset h\circ W.
	\end{align*}. 
	Since $A$ is maximal $h$ symmetric, we have $(x,0)\in A$ and $x\in\ker A$.
	Then we have $\ker A\supset(\image A)^{\Omega,l}$. By Lemma \ref{l:ker-dom-symmetric}, we have $\ker A\subset(\image A)^{\Omega,l}$ and then $\ker A=(\image A)^{\Omega,l}$. Apply the result for $A^{-1}$, we have $A0=(\dom(A))^{\Omega,r}$.
	\newline 2. Let $(x,y)$ be in $A^\Omega$. Then we have \begin{align*}
		0=\omega((x,y),(0,y_1))=\Omega(x,y_1)
	\end{align*}
    for all $y_1\in A0$. By Step 1, we have $x\in (A0)^{\Omega,l}=((\dom(A))^{\Omega,r})^{\Omega,l}.$ Then we obtain    
    $\dom(A^\Omega)\subset((\dom(A))^{\Omega,r})^{\Omega,l}$. Apply the result for $A^{-1}$, we have $\image A^\Omega\subset((\image A)^{\Omega,l})^{\Omega,r}$.    
    \newline 3. Let $(x,y)$ be in $A^\Omega$ with $x\in\dom(A)$. Then there is a $y_2\in Y$ such that $(x,y_2)\in A$. Since $h\circ A\subset A^\Omega$, there holds $(0,y-hy_2)=(x,y)-(x,hy_2)\in A^\Omega$. For each $(x_3,y_3)\in A$ we have 
    \begin{align*}
    	0=\omega((x_3,hy_3),(0,y-hy_2))=\Omega(x_3,y-hy_2).
    \end{align*}
    By Step 1, we have $y-hy_2\in  (\dom(A))^{\Omega,r}=A0$ and $y\in hy_2+A0=hAx$. Then we have $(x,y)\in h\circ A$ and $A^\Omega\cap(\dom(A)\times Y)\subset h\circ A$. Note that $A^\Omega\cap(\dom(A)\times Y)\supset hA$. Then we get 
    $A^\Omega\cap(\dom(A)\times Y)=hA$ and $A^\Omega0=A0$. Apply the result for $A^{-1}$, we have $\ker A^\Omega=\ker A$. 
    \newline 4. Since $|h|=1$ and $A$ is $h$ symmetric with respect to $\Omega$, we have $\ker A\subset\ker A$. By Step 1, we have $\ker A=(\ran A)^{\Omega,l}$. 
    Since $\ker A\subset\ker Q_A\subset (\ran A)^{\Omega,l}$, we have $\ker A=\ker Q_A$.    
\end{proof}

\begin{corollary}\label{e:max-symmetric-selfadjoint}
	Let $X$ and $Y$ be two vector space over $\KK$. Let $\Omega\colon X\times Y\to\KK$ be a nondegenerate sesquilinear form. Let $h\in\KK$ be such that $|h|=1$. 
	Let $A\colon X\leadsto Y$ be a maximal $h$ symmetric linear relation with respect to $\Omega$. Then $A$ is $h$ selfadjoint if one of the following four conditions holds:
	\begin{itemize}
		\item [(i)] $\ran A=\ran A^\Omega$,
		\item[(ii)] $\ran A=((\image A)^{\Omega,l})^{\Omega,r}$,
		\item [(iii)] $\dom(A)=\ran (A)^\Omega$,
		\item[(iv)] $\dom(A)=((\dom(A))^{\Omega,r})^{\Omega,l}$.
	\end{itemize}  
\end{corollary}

\begin{proof}
	By \eqref{e:im-dom-1} and $\ran A\subset\ran A^\Omega$, we have (ii)$\Rightarrow$(i).
	
	Assume (i) holds. By \eqref{e:ker-dom}, we have $\ker A=\ker A^\Omega$. Since $h\circ A\subset A^\Omega$ and $|h|=1$, we have 
	\begin{align*}
		h\circ A=A^\Omega\cap(X\times \image A)=A^\Omega.
	\end{align*}
	
	Apply the arguments to $A^{-1}$, we obtain the rest two implications.
\end{proof}

The following example gives a symmetric operator $A$ such that $\ker A\subsetneq\ker Q^A$.

\begin{example}\label{ex:no-max-symmetric}
	We define a linear operator $A:\R^2\supset\R\times\{0\}\to\R^2$ by
	$A(x,0)=(0,x)$ for each $x\in\R$. Then $A$ is symmetric, and $\ker A=\{0\}\subsetneq\R\times\{0\}=\ker Q_A$.
\end{example}

The index of a $h$-symmetric linear relation ($|h|=1$) is always non-positive. More precisely, we have the following results.

\begin{lemma}\label{l:selfadjoint-index-non-positive}
	Let $X$ and $Y$ be two vector space over $\KK$. Let $\Omega\colon X\times Y\to\C$ be a nondegenerate sesquilinear form. Let $h\in\KK$ be such that $|h|=1$. Let $A\colon X\leadsto Y$ be an $h$ symmetric linear relation with respect to $\Omega$. Then the following hold.
	\newline (a) We have $\Index A\le 0$ and $\dim\ker A<+\infty$ if one of $\dim\ker A$ and $\dim Y/\image A$ is finite.
	\newline (b) The equality $\Index A=0$ holds if and only if $\dim\ker A=\dim\ker A^{\Omega}$ and $\image A=((\image A)^{\Omega,l})^{\Omega,r}$.
	\newline (c) The relation $A$ is $h$ selfadjoint if the equality $\Index A=0$ holds.
\end{lemma}

\begin{proof} Let $\omega$ be defined by \eqref{e:natural-symplectic-structure}.%
	\newline (a) Firstly assume that $\KK=\C$. Then $a\circ A$ is symmetric for $a\in\C$ with $a^2=h$. Note that $a\circ A+X\times\{0\}=X\times\image A$. If one of $\dim\ker A$ and $\dim Y/\image A$ is finite,
	by \cite[Lemma 1.2.8]{BoZh18} we have
	\begin{align*}
		\dim\ker A&= \dim(A\cap(X\times\{0\}))\\
		&\le\dim((a\circ A)^\omega\cap(X\times\{0\})^\omega)\\
		&\le\dim (X\times Y)/((a\circ A)+X\times\{0\})\\
		&=\dim Y/\image A.
	\end{align*}
    Then we have $\Index A\le 0$ and $\dim\ker A<+\infty$.
    
    For $\KK=\R$, we apply the complex case on $A\otimes_{\R}\C$ and get the results.
    \newline (b) By \cite[Lemma 1.2.8]{BoZh18}, $\Index A=0$ holds if and only if  $\dim\ker A=\dim\ker A^\omega$ and $(X\times\image A)^{\omega\omega}=X\times\image A$, if and only if $\dim\ker A=\dim\ker A^\omega$ and $\image A=(\image A^{\Omega,l})^{\Omega,r}$.
    \newline (c) For $\KK=\C$, we note that $a\circ A$ is symmetric for $a\in\C$ with $a^2=h$. and $\Index(a\circ A)=0$. By \cite[Proposition 1 ]{BoZh13}, $a\circ A$ is selfadjoint. Then $A$ is $h$-selfadjoint. 
    
    For $\KK=\R$, by the complex case the complexified linear relation $a\circ A\otimes_{\R}\C$ is selfadjoint. Then $A$ is $h$  selfadjoint. 
\end{proof}
	
We have the following generalization of \cite[Lemma 2.2]{CCGP23}, where the full proof of the corresponding result in N.J. Kalton, R. C. Swanson \cite{KaSW82} was given. 

\begin{proposition}\label{p:index-selfadjoint}
	Let $X$ be a Banach space over $\KK$. Let $h\in\KK$ be such that $|h|=1$. Let $A\colon X\leadsto X^*$ be an $h$ selfadjoint linear relation. Then the following hold. 
	\begin{itemize}
		\item[(a)] Assume that $A$ is semi-Fredholm. Then $\Index A=0$ holds if and only if $\image\iota_X\supset\dom(A^*)$, if and only if $\image A=\image A^*$. In particular, $\Index A=0$ holds if $X$ is reflexive.
		\item[(b)] We have $A^*0=A0$, and $\iota_X$ induces a linear isomorphism $X/\dom(A)\to \image\iota/(\image\iota\cap\dom(A^*))$. 
	\end{itemize}
\end{proposition}

\begin{proof}  
	1. Since $|h|=1$, $h\circ A=A^\omega=A^*\circ\iota_X$ and $\iota_X$ is an injection, by \cite[Exercise B.2.11]{GWh78} for $G=\iota_X$ and $F=A^*$, 
	there are exact sequences
	\begin{align}\label{e:exact-sequence-A}	
		\begin{split}
			\{0\}&\to\ker A\to\ker A^*
			\to \dom(A^*)/(\dom(A^*)\cap\image \iota_X)\\
			&\to X^*/\image A\to X^*/\image A^*\to \{0\},
		\end{split}			
	\end{align}
	\begin{align}\label{e:exact-sequence-B}	
		\begin{split}
			\{0\}&\to A^*0\to A0
			\to \{0\}\\
			&\to \image\iota/(\image\iota\cap\dom(A^*))\to X/\dom(A)\to \{0\},
		\end{split}			
	\end{align}
 	\newline 2. {\em Proof of (a).}
 	
 	Since $A$ is semi-Fredholm,
	by Lemma \ref{l:selfadjoint-index-non-positive} and the proof of \cite[Theorem IV.5.13]{Ka95}, we have 
	\begin{align*}
		\Index A&=\dim\ker A-\dim\ker A^*\\
		&=\dim X^*/\image A^*-\dim X^*/\image A\le 0.
	\end{align*}
	By \eqref{e:exact-sequence-A}, $\Index A=0$ holds if and only if the injection $\ker A\to\ker A^*$ is an isomorphism, if and only if $\dom(A^*)\cap\image \iota_X=\dom(A^*)$, if and only if $\image\iota_X\supset\dom(A^*)$, if and only if the surjection $X^*/\image A\to X^*/\image A^*$ is an isomorphism, if and only if $\image A=\image A^*$.	
	\newline 3. By \eqref{e:exact-sequence-A}, we obtain (b).	
\end{proof}

For symplectic linear subspace, we have the following.

\begin{lemma}\label{l:symplectic-sum}
	Let $(X,\omega)$ be a symplectic vector space with a symplectic linear subspace $X_0$.
	Then $\image J(\cdot)|_{X_0}\subset \image J_{\omega|_{X_0\times X_0}}$ holds if and only if
	\begin{align*}
		X=X_0\oplus X_0^{\omega}.
	\end{align*}
    In this case, we have $X_0^{\omega\omega}=X_0$.
\end{lemma}

\begin{proof}
	Since $X_0$ is symplectic,  $J_{\omega|_{X_0\times X_0}}$ is injective and we have $X_0\cap X_0^{\omega}=\{0\}$. Then we have
	\begin{align*}
		&\image J(\cdot)|_{X_0}\subset \image J_{\omega|_{X_0\times X_0}}\\
		\Leftrightarrow&\forall x\in X,\exists x_0\in X_0\text{ such that }J(x)|_{X_0}=J_{\omega|_{X_0\times X_0}}(x_0)\\
		\Leftrightarrow&\forall x\in X,\exists x_0\in X_0\text{ such that }x-x_0\in X_0^{\omega}\\
		\Leftrightarrow&X\subset X_0+X_0^{\omega}\\
		\Leftrightarrow&X= X_0\oplus X_0^{\omega}.
	\end{align*}
    In this case, we have $X_0^{\omega}\cap X_0^{\omega\omega} =(X_0+X_0^\omega)^\omega=\{0\}$.
    By \cite[Lemma A.1.1]{BoZh18}, we have
    \begin{align*}
    	X_0^{\omega\omega}
    	&=X_0^{\omega\omega}\cap(X_0+X_0^{\omega})\\
    	&=X_0+X_0^{\omega\omega}\cap X_0^{\omega}=X_0.    	
    \end{align*}
\end{proof}

\begin{corollary}\label{c:symplectic-sum} (cf. \cite[Lemma 2.7]{CCGP23})
	Let $(X,\omega)$ be a symplectic Banach space with a closed symmplectic linear subspace $X_0$. Then the following hold.
	\newline (a) If $(X_0,\omega|_{X_0\times X_0})$ is a strong symplectic Banach space, we have
	\begin{align*}
		X=X_0\oplus X_0^{\omega}.
	\end{align*}
    \newline (b) Assume that $(X,\omega)$ is a strong symplectic Banach space and $X_0+X_0^\omega$ is closed. Then $(X_0,\omega|_{X_0\times X_0})$ is a strong symplectic Banach space.
\end{corollary}

\begin{proof}
	(a) Since $(X_0,\omega|_{X_0\times X_0})$ is a strong symplectic Banach space, we have $\image J_{\omega|_{X_0\times X_0}}=X_0^*\supset \image J(\cdot)|_{X_0}$. By Lemma \ref{l:symplectic-sum}, we have $X=X_0\oplus X_0^\omega$.
	\newline (b) By \cite[Lemma 1.2.6]{BoZh18}, we have 
	\begin{align*}
		X_0+X_0^\omega=(X_0+X_0^\omega)^{\omega\omega}=(X_0^\omega\cap X_0^{\omega\omega})^\omega=(X_0^\omega\cap X_0)^\omega=\{0\}^\omega=X.
	\end{align*}
	
	Since $(X,\omega)$ is a strong symplectic Banach space, we have $\image J=X^*$. By Hahn-Banach theorem, we have $X_0^*\subset\image J(\cdot)|_{X_0}$. 	
	Since $X=X_0\oplus X_0^{\omega}$ holds, by Lemma \ref{l:symplectic-sum}, we have $\image J(\cdot)|_{X_0}\subset\image J_{\omega|_{X_0\times X_0}}$. Since $\image J(\cdot)|_{X_0}$ and $\image J_{\omega|_{X_0\times X_0}}$ are subsets of $X_0^*$, we have $X_0^*=\image J(\cdot)|_{X_0}=\image J_{\omega|_{X_0\times X_0}}$, i.e. $X_0$ is a strong symplectic Banach space.
\end{proof}

\subsection{Perturbed Morse index}\label{ss:stable-morse}

Before study the stability theorems of maximal isotropic spaces, we study the perturbed Morse index for symmetric forms with varying domains in a fixed Banach space. Note that in general one can not make the domains fixed.

Let $V$ be a vector space over $\KK$ and $Q\colon V\times V\to\KK$ be a symmetric form.
For each subset $\lambda$ of $V$, we denote by 
\begin{align*}
	\lambda^Q\ :=\lambda^{Q,l}=\lambda^{Q,r}.
\end{align*}
Then $m^{\pm}(Q)$ and $m^0(Q)$ denotes the {\em Morse positive (or negative) index} and the {\em nullity} of $Q$ respectively.  

\begin{definition}\label{d:bounded-symmetric-pair}
	(a) Let $X$ be a Banach space over $\KK$ with a closed linear subspace $V$. Let $Q\colon V\times V\to\KK$ be a symmetric form. If 
	\begin{align}\label{e:norm-Q}
		\|Q\|\ :=\begin{cases}
			\sup\limits_{x,y\in V\setminus\{0\}}\frac{|Q(x,y)|}{\|x\|\|y\|}<+\infty&\text{ if } V\ne\{0\},\\0&\text{ if } V=\{0\},
		\end{cases}
	\end{align}
	we call $(Q,V)$ a {\em bounded symmetric pair} of $X$. 
	\newline (b) Let $X$ be a Banach space over $\KK$ with a closed linear subspace $V$. Let $Q\colon V\times V\to\KK$ be a bounded semi-definite symmetric form. The {\em reduced minimum modulus} $\gamma(Q)$ of $Q$ is defined by 
	\begin{align}\label{e:reduced-minimum-modulus-Q}
		\gamma(Q)\ :=\begin{cases}
			\inf\limits_{x\in V\setminus V^Q}\frac{|Q(x,x)|}{(\dist(x,V^Q))^2}&\text{ if } V\ne V^Q,
			\\0&\text{ if } V=V^Q,
		\end{cases}
	\end{align}
	\newline (c) Let $X$ be a Banach space with two  symmetric pairs $(Q,V)$ and $(R,W)$. Let $c\ge 0$ be a real number. We define the {\em $c$-gap} $\delta_c(Q,R)$ between $R$ and $Q$ to be the infimum of the non-negative number $\delta$ such that 
	\begin{align}\label{e:c-gap-pair}
		\begin{split}
			|Q(x,y)-&R(u,v)|\le \delta(\|u\|+\|x\|)(\|v\|+\|y\|)+\\
			& c\left((\|u\|+\|x\|)\|v-y\|+\|u-x\|(\|v\|+\|y\|)\right)
		\end{split}    	
	\end{align}
	for all $x,y\in V$ and $u,v\in W$.
\end{definition}

We have the following properties of the $c$-gap.

\begin{lemma}\label{l:c-gap}
	Let $X$ be a Banach space with two symmetric pairs $(Q,V)$ and $(R,W)$. Let $c\ge 0$ be a real number. Then the following hold.
	\begin{itemize}
		\item[(a)] The inequality (\ref{e:c-gap-pair}) holds for $\delta=\delta_c(Q,R)$ and $\delta_c(Q,R)=\delta_c(R,Q)$.
		\item[(b)] Assume that $\delta_c(Q,R)$ is finite and $V,W$ are closed. Then $(Q,V)$ and $(R,W)$ are bounded and we have 
		\begin{align}\label{e:norm-Q-gap-pair}
			\max\{\|Q\|,\|R\|\}\le\delta_c(Q,R)+2c.
		\end{align}
		\item[(c)] Assume that $(Q,V)$ and $(R,W)$ are bounded. Then we have 
		\begin{align}\label{e:range-delta-c}
			0\le \delta_c(Q,R)\le\max\{\|Q\|,\|R\|\}.
		\end{align}
		\item[(d)] The function $\delta_c(Q,R)$ is  decreasing on $c\ge 0$.
		\item[(e)] Let $M$ and $N$ be linear subspace of $V$ and $W$ respectively. Then we have 
		\begin{align}\label{e:monotone-delta-c}
			\delta_c(Q|_M,R|_N)\le\delta_c(Q,R).
		\end{align}
		\item[(f)] Assume that $(Q,V)$ and $(R,W)$ are bounded. If $V=W$ holds, we have 
		\begin{align}\label{e:c-gap-whole-domain}
			\frac{\|R-Q\|}{4}\le\delta_c(Q,R)\le\|R-Q\|
		\end{align}
		for $c\ge\|Q\|$.
	\end{itemize}
\end{lemma}

\begin{proof} (a) 
	By definition, we have $\delta_c(Q,R)=\delta_c(R,Q)$.
	
	Let $\varepsilon>0$ be arbitrary. For $x,y\in V$ and $u,v\in W$, \eqref{e:c-gap-pair} holds for $\delta=\delta_c(Q,R)+\varepsilon$. On letting $\varepsilon\to 0$, we have \eqref{e:c-gap-pair} holds for $x,y\in V$, $u,v\in W$ and $\delta=\delta_c(Q,R)$.
	\newline
	(b) Take $u=v=0$ in \eqref{e:c-gap-pair}, we have
	\begin{align*}
		|Q(x,y)|\le(\delta_c(Q,R)+2c)\|x\|\|y\|
	\end{align*}
	for all $x,y\in V$.
	By definition, we have  $\|Q\|\le\delta_c(Q,R)+2c$. 
	By changing the role of $Q$ and $R$, we have $\|R\|\le\delta_c(Q,R)+2c$.
	\newline
	(c) By definition, we have $\delta_c(Q,R)=\delta_c(R,Q)\ge 0$.
	Let $x,y\in V$ and $u,v\in W$ be elements in $X$. Then we have
	\begin{align*}
		|Q(x,y)&-R(u,v)|\le|Q(x,y)|+|R(u,v)|\\
		&\le \|Q\|\|x\|\|y\|+\|R\|\|u\|v\|\\
		&\le\max\{\|Q\|,\|R\|\}(\|u\|+\|x\|)(\|v\|+\|y\|).		
	\end{align*}
	By definition we obtain \eqref{e:range-delta-c}.
	\newline (d) Let $c$, $d$ be non-negative numbers with $c\le d$. Then we have
	\begin{align*}
		|Q(x,y)-&R(u,v)|\le \delta_c(P,Q)(\|u\|+\|x\|)(\|v\|+\|y\|)+\\
		& d\left((\|u\|+\|x\|)\|v-y\|+\|u-x\|(\|v\|+\|y\|)\right)
	\end{align*}
	for all $x,y\in V$ and $u,v\in W$. By definition, we have  $\delta_d(P,Q)\le\delta_c(P,Q)$. 
	\newline (e) Since \eqref{e:c-gap-pair} holds for
	all $x,y\in M$ and $u,v\in N$, by definition, we have $\delta_c(Q|_M,R|_N)\le\delta_c(P,Q)$.
	\newline (f) Let $x$, $y$, $u$, $v$ be in $V=W$.
	Since $c\ge\|Q\|$ holds, we have 
	\begin{align*}
		|Q(x,y)-&R(u,v)|\le |Q(u,v)-R(u,v)|+|Q(x,y-v)|+|Q(x-u,v)|\\
		\le&\|Q-R\|\|u\|\|v\|+\|Q\|(\|x\|\|v-y\|+\|u-x\|\|v\|)\\
		\le&
		\|Q-R\|(\|u\|+\|x\|)(\|v\|+\|y\|)+\\
		& c\left((\|u\|+\|x\|)\|v-y\|+\|u-x\|(\|v\|+\|y\|)\right).
	\end{align*}
	By definition, we have 
	\begin{align*}
		\delta_c(Q,R)\le\|R-Q\|.
	\end{align*}
	By \eqref{e:c-gap-pair}, we have
	\begin{align*}
		|Q(x,y)-R(x,y)|\le 4(\delta_c(Q,R))\|x\|\|y\|
	\end{align*}
	for all $x,y\in V=W$. Then we have
	\begin{align*}
		\delta_c(Q,R)\ge\frac{\|R-Q\|}{4}.
	\end{align*}
\end{proof}

The main result of this subsection is the following.

\begin{proposition}\label{p:perturbed-Morse-index}
	Let $X$ be a Banach space with two bounded symmetric pairs $(Q,V)$ and $(R,W)$. Let $c\ge 0$ be a real number.
	Let $\alpha$ be a linear subspace of $V$ with $\dim\alpha=k\in[1,+\infty)$, $h$ be in $\{1,-1\}$ and $h Q|_\alpha$ be positive definite. Assume that 
	\begin{align}\label{e:perturbed-Morse-index-delta}
		k(\delta_c(Q,R)+
		2c\delta(V,W))(2+\delta(V,W))
		<\gamma(Q|_\alpha).
	\end{align}
	Then we have
	\begin{align}\label{e:perturbed-Morse-index}
		m^+(hR)\ge \dim\alpha.
	\end{align} 	    
\end{proposition}

\begin{proof}
	Take an orthogonal bases $v_1,\ldots,v_k$ of $\alpha$ with respect to $h Q$ such that $\|v_i\|=1$ for each $i=1,\ldots,k$. For each $\varepsilon>0$ and $i=1,\ldots, k$,  there is a $w_i\in W$ with $\|v_i-w_i\|\le(1+\varepsilon)\delta(V,W)$. 
	
	By \eqref{e:c-gap-pair}, we have
	\begin{align*}
		|R(w_i,w_j)&-Q(v_i,v_j)|\le \delta_c(Q,R)(\|w_i\|+\|v_i\|)(\|w_j\|+\|v_j\|)+\\
		&c\left(\|v_i-w_i\|(\|v_j\|+\|w_j\|)+\|v_j-w_j\|(\|v_i\|+\|w_i\|)\right)\\
		\le&\delta_c(Q,R)(2+(1+\varepsilon)\delta(V,W))^2
		+
		\\&2c(1+\varepsilon)\delta(V,W)(2+(1+\varepsilon)\delta(V,W)):\ =C(\varepsilon).
	\end{align*}
	Set $A\ :=(Q(v_i,v_j))_{1\le i,j\le k}$ and $B\ :=(R(w_i,w_j))_{1\le i,j\le k}$. Then we have $hA\ge\gamma(Q|_\alpha)I_k>0$ and $\|B\|\le kC(\varepsilon)$ (here we view $B$ as a sesquilinear form on $\KK^k$). Since $kC(0)<\gamma(Q|_\alpha)$, for $\varepsilon>0$ sufficiently small, we have \begin{align*}
		h(B-A)\ge -kC(\varepsilon)I_k>-\gamma(Q|_\alpha)I_k.
	\end{align*} 
	Since $h A\ge \gamma(Q|_\alpha)I_k>0$, we have $hB>0$. Then for each $(a_1,\ldots,a_k)\in\KK^k\setminus\{0\}$ and $w\ :=\sum_{i=1}^ka_iw_i$, we have $hR(w,w)>0$. Therefore $\{w_1,\ldots,w_k\}$ is a linearly independent subset and we have $m^+(hR)\ge k$.    
\end{proof}

\subsection{Stability theorems of maximal isotropic subspaces}\label{ss:stable-max-isotropic}

In this subsection we study the stability theorems of maximal isotropic subspaces. 
Firstly we recall the concept of symplectic reduction.

Let $(X,\omega)$ be a symplectic vector space with an isotropic subspace $\lambda$. Then the induced skew symmetric form on the reduced space $\lambda^\omega/\lambda$ is defined by
\begin{align}\label{e:red-structure}
	\wt{\omega}(x+\lambda,y+\lambda)\ :=\ \omega(x,y) \text{ for all } x,y\in \lambda^\omega.
\end{align}
The space $(\lambda^\omega/\lambda, \wt{\omega})$ is symplectic if and only if $\lambda=\lambda^{\omega\omega}$. We call
$(\lambda^\omega/\lambda, \wt{\omega})$ the {\em
symplectic reduction} of $X$ by $\lambda$ and denote by $\pi_\lambda\colon\lambda^\omega\to \lambda^\omega/\lambda$ the canonical projection.
Let $\alpha$ be a linear subspace of $X$.
The {\em symplectic reduction} of $\alpha$ by $\lambda$ is defined by
\begin{align}\label{e:red-subspace}
	\pi_\lambda(\alpha)\ :=\ \bigl((\alpha+\lambda)\cap \lambda^\omega\bigr)/\lambda\ =\ \bigl(\alpha\cap  \lambda^\omega+\lambda\bigr)/\lambda.
\end{align}

We need the following properties of maximal isotropic subspaces.

\begin{proposition}\label{p:maximal-isotropic-properties}
	Let $(X,\omega)$ be a symplectic vector space with an isotropic subspace $\lambda$. Then the following hold.
	\newline (a) If $\lambda$ is a maximal isotropic subspace of $(X,\omega)$, we have $\lambda=\lambda^{\omega\omega}$.
	\newline (b) $\lambda$ is a maximal isotropic subspace of $(X,\omega)$ if and only if 
	\begin{align}\label{e:maximal-isotropic}
		\lambda=\{x\in\lambda^\omega;\;\omega(x,x)=0\}.
	\end{align} 		
	\newline (c) Assume that $(X,\omega)$ is a complex symplectic vector space. Then $\lambda$ is a maximal isotropic subspace of $(X,\omega)$ if and only if $\mi\tilde\omega$ is definite on $\lambda^\omega/\lambda$.
	\newline (d) Assume that $(X,\omega)$ is a real symplectic vector space. Then $\lambda$ is a maximal isotropic subspace of $(X,\omega)$ if and only if $\lambda=\lambda^\omega$.	
\end{proposition}

\begin{proof} 
	(a) Since $\lambda\subset\lambda^\omega$ holds, we have 
	\begin{align*}
		\lambda\subset\lambda^{\omega\omega}\subset\lambda^\omega=\lambda^{\omega\omega\omega}.
	\end{align*} 
    Since $\lambda$ is a maximal isotropic subspace of $(X,\omega)$, we have $\lambda=\lambda^{\omega\omega}$.
    \newline (b) (i) {\em Necessity.} 
    
    Since $\lambda$ is isotropic, we have $\lambda\subset\{x\in\lambda^\omega;\;\omega(x,x)=0\}$. Let $x\in\lambda^\omega$ be such that $\omega(x,x)=0$.     
    Assume that $x\notin\lambda$. Since $\omega(x,x)=0$, the subspace $\lambda\oplus\Span\{x\}\supsetneq\lambda$ is an isotropic subspace of $(X,\omega)$. This contradicts to the fact that $\lambda$ is a
    maximal isotropic subspace of $(X,\omega)$.
    Then we have $\lambda\supset\{x\in\lambda^\omega;\;\omega(x,x)=0\}$. Thus we obtain \eqref{e:maximal-isotropic}.
    
    (ii) {\em Sufficiency.}
    
    Let $\mu\supset\lambda$ be isotropic. Then we have $\mu\subset\mu^\omega\subset\lambda^\omega$.
    By \eqref{e:maximal-isotropic}, we have $\mu\subset\{x\in\lambda^\omega;\;\omega(x,x)=0\}=\lambda$. Then we have $\lambda=\mu$. Since $\mu$ is a arbitrary isotropic subspace containing $\lambda$, $\lambda$ is maximal isotropic.
    \newline (c) (i) {\em Necessity.} 
    
    By (a), $\tilde\omega$ is a symplectic form on $\lambda^\omega/\lambda$. If $\mi\tilde\omega$ is not definite on $\lambda^\omega/\lambda$, by (b), there are $x,y\in\lambda^\omega$ such that $\mi\omega(x,x)>0$ and $\mi\omega(y,y)<0$.
    Then there is an $a\in\C\setminus\{0\}$ such that $\omega(ax,y)=a\omega(x,y)\in\R$ and
    \begin{align*}
    	\mi\omega(ax+y,ax+y)=|a|^2\mi\omega(x,x)+\mi\omega(y,y)=0.
    \end{align*} 
    By (b), we have $ax+y\in\lambda$. Then we have \begin{align*}
    	\omega(y,y)=\omega(ax+y-ax,ax+y-ax)=|a|^2\omega(x,x).
    \end{align*} 
    From this we obtain that $|a|^2\mi\omega(x,x)=\mi\omega(y,y)$, which contradicts to the facts $\mi\omega(x,x)>0$ and $\mi\omega(y,y)<0$.
    
    (ii) {\em Sufficiency.} 
    
    Since $\mi\tilde\omega$ is definite on $\lambda^\omega/\lambda$, \eqref{e:maximal-isotropic} holds. By (b), $\lambda$ is maximal isotropic.
    \newline (d) (i) Necessity.
    
    Since $X$ is real, we have $\omega(x,x)=0$ for all $x\in X$. By (b), $\lambda$ is a Lagrangian subspace of $(X,\omega)$.
    
    (ii) Sufficiency. 
    
    Since $\lambda$ is a Lagrangian subspace of $(X,\omega)$, it is maximal isotropic.     
\end{proof}

Based on Proposition \ref{p:maximal-isotropic-properties}, we can define the sign of a maximal isotropic subspace.

\begin{definition}\label{d:sign-max-isotroic}
	Let $(X,\omega)$ be a symplectic vector space and $\lambda$ be a maximal isotropic subspace of $(X,\omega)$. We define the {\em sign} $h_\lambda$ of $\lambda$ by
	\begin{align}\label{e:sign-max-isotroic}
		h_\lambda\ :=
		\begin{cases}
			\pm 1&\text{ if }\pm\mi\omega|_\lambda\text{ is positive definite},\\
			0&\text{ if }\lambda=\lambda^\omega.
		\end{cases}
	\end{align}
	If $(X,\omega)$ is a symplectic Banach space, by Proposition \ref{p:maximal-isotropic-properties}.a, the space $\lambda=\lambda^{\omega\omega}$ is closed. We define the {\em reduced minimum modulus} $\gamma_\lambda$ of $\lambda$ by 
	\begin{align}\label{e:delta-lambda}
		\gamma_\lambda\ :=
		\begin{cases}
			\inf\limits_{x\in\lambda^\omega\setminus\lambda}\frac{\mi h_\lambda\omega(x,x)}{(\dist(x,\lambda))^2} &\text{ if }\lambda\subsetneq\lambda^\omega\\
			0&\text{ if } \lambda=\lambda^\omega
		\end{cases}
	\end{align}
	respectively.
\end{definition} 

By definition, a maximal isotropic subspace $\lambda$ is Lagrangian if and only if $h_\lambda=0$.

The starting point of our stability theorems of maximal isotropic subspaces is the following result of P. Hess and T. Kato \cite[Lemma 1]{HeKa70}).

\begin{lemma}\label{l:relative-dimension-0}
	Let $X$ be a Banach space with three closed linear subspaces $M$, $N$ and $N'$, where $N'\subset N$. If there hold that 
	\begin{align*}
		(1+\delta(N,M))	(1+\delta(M,N'))<2,
	\end{align*}
	we have $N'=N$.
\end{lemma}

\begin{proof}[Proof of Theorem \ref{t:stable-max-isotropic}] 		
	Assume that $\mu$ is not a a maximal isotropic subspace of $(X,\omega)$ with $h_\mu\in\{h_\lambda,0\}$. Then we have
	$\mu\subsetneq\mu^\omega$. 	
	By Proposition \ref{p:maximal-isotropic-properties}.b, there is an $x_1\in\mu^\omega\setminus\mu$ such that 
	\begin{align*}
		\mi h_\lambda\omega(x_1,x_1)\le 0.
	\end{align*}
	Since $\mu$ is closed, by \cite[Lemma III.1.12]{Ka95}, for each $\varepsilon\in(0,1)$, there is a $y_1\in\mu$ such that 
	\begin{align*}
		\dist(x_1,\mu)>(1-\varepsilon)\|x_1-y_1\|>0.
	\end{align*}
	Set $x\ :=\frac{x_1-y_1}{\|x_1-y_1\|}$. Then we have $x\in\mu^\omega\setminus\mu$, $\|x\|=1$, $\mi h_\lambda\omega(x,x)\le 0$, and
	\begin{align*}
		\dist(x,\mu)=\dist(\frac{x_1}{\|x_1-y_1\|},\mu)>1-\varepsilon>0.
	\end{align*}
	
	By definition, there is a $y\in\lambda^{\omega_0}$ such that
	\begin{align*}
		\|x-y\|<\delta(\mu^\omega,\lambda^{\omega_0})+\varepsilon.
	\end{align*}    
	If $h_\lambda=0$, we have $\lambda=\lambda^{\omega_0}$ and $\dist(x,\lambda)=0$. If $h_\lambda\ne0$, we have $\gamma_\lambda>0$. By \eqref{e:delta-lambda}, we have 
	\begin{align*}
		\dist(y,\lambda)\le&\gamma_\lambda^{-\frac{1}{2}}(\mi h_\lambda\omega_0(y,y))^{\frac{1}{2}}\\
		=&\left(\mi h_\lambda\gamma_\lambda^{-1}\left(\omega(x,x)+\omega_0(x,x)-\omega(x,x)+\right.\right.\\
		&\left.\left.\omega_0(y-x,y)+\omega_0(x,y-x)\right)\right)^{\frac{1}{2}}\\
		\le& \left(|h_\lambda|\gamma_\lambda^{-1}\left(\|L_\omega-L_{\omega_0}\|+\|L_{\omega_0}\|\|x-y\|(2+\|x-y\|)\right)\right)^{\frac{1}{2}}.
	\end{align*}
	By \cite[Lemma IV.2.2]{Ka95}, we have
	\begin{align*}
		(1+\delta(\lambda,\mu))&(\|x-y\|+\dist(y,\lambda))
		\ge(1+\delta(\lambda,\mu))\dist(x,\lambda)\\
		&\ge\dist(x,\mu)-\delta(\lambda,\mu)\\
		&>1-\varepsilon-\delta(\lambda,\mu).
	\end{align*}
	Then we have
	\begin{align*}
		\delta(\mu^\omega,\lambda^{\omega_0})+&\varepsilon+\left(|h_\lambda|\gamma_\lambda^{-1}\left(\|L_\omega-L_{\omega_0}\|+\|L_{\omega_0}\|(\delta(\mu^\omega,\lambda^{\omega_0})+\varepsilon)\right.\right.\\
		&\left.\left.(2+\delta(\mu^\omega,\lambda^{\omega_0})+\varepsilon)\right)\right)^{\frac{1}{2}}
		>\frac{1-\varepsilon-\delta(\lambda,\mu)}{1+\delta(\lambda,\mu)}.
	\end{align*} 
	On letting $\varepsilon\to 0$ we obtain
	\begin{align}\label{e:non-maximal-isotropic}   
		\delta(\mu^\omega,\lambda^{\omega_0})+l(\mu,\omega)
		\ge\frac{1-\delta(\lambda,\mu)}{1+\delta(\lambda,\mu)}.
	\end{align} 
	This contradicts to \eqref{e:stable-maximal-isotropic}. So $\mu$ is a maximal isotropic subspace of $(X,\omega)$ with $h_\mu\in\{h_\lambda,0\}$.
\end{proof}

\begin{proof}[Proof of Theorem \ref{t:family-max-isotropic}]
	1. Denote by 
	\begin{align*}		
		K_1\ :=&\{s\in J;\;\lambda(s)=\lambda^{\omega(s)\omega(s)} \},\\
		K_2\ :=&\{s\in J;\;\lambda(s)=\lambda^{\omega(s)} \}.
	\end{align*}
	Then $K_1$ and $K_2$ are closed subsets of $J$. By Proposition \ref{p:maximal-isotropic-properties}.a, we have $s_0\in K_1$. By Lemma \ref{l:relative-dimension-0}, $K_1$ and $K_2$ are open subsets of $J$. Since $J$ is connected, we have $K_1=J$ and $K_2$ is either $J$ or empty set. Then we have $K_2=J$ if $s_0\in K_2$, and in this case we have $h_{\lambda(s)}=h_{\lambda(s_0)}=0$ holds for each $s\in J$. 
	\newline 2. Assume that $K_2=\emptyset $ holds. Then we have $h_{\lambda(s_0)}\ne 0$. By Proposition \ref{p:maximal-isotropic-properties}.d, we have $\KK=\CC$. Denote by 
	\begin{align*}
		K_3\ :=\{s\in J;\;m^-(\mi h_{\lambda(s_0)}\omega_s|_{\lambda(s)^\omega(s)})=0 \}.
	\end{align*}
	By Proposition \ref{p:perturbed-Morse-index}, $J\setminus K_3$ is open and $K_3$ is closed.
	Since $K_1=J$, by Proposition \ref{p:maximal-isotropic-properties}.c, 
	we have \begin{align*}
		K_3=&\{s\in J;\;\lambda(s)\text{ is a maximal isotropic subspace of }\\
		&(X,\omega(s))\text{ and } h_{\lambda(s)}\in\{h_{\lambda(s_0)},0\}\}\\
		=&\{s\in J;\;\lambda(s)\text{ is a maximal isotropic subspace of }\\
		&(X,\omega(s))\text{ and } h_{\lambda(s)}=h_{\lambda(s_0)}\}.
	\end{align*}
	Then we have $s_0\in K_3$. By Theorem \ref{t:stable-max-isotropic}, $K_3$ is open. Since $K_3\subset J$ is a nonempty set and $J$ is connected, we have $K_3=J$. 
\end{proof}

\subsection{Stability theorems in strong symplectic Banach spaces}\label{ss:stable-max-isotropic-strong}
In strong symplectic Banach spaces, some conditions in stability theorems of maximal isotropic subspaces hold. Then these stability theorems can be simplified.

If $(X,\omega(s))$ is a strong symplectic Banach space for each $s\in J$, the continuity of the family $\{\lambda(s)^{\omega(s)}\}_{s\in J}$ follows from that of the family $\{\lambda(s)\}_{s\in J}$ and
the following estimate (cf. \cite[Theorem IV.2.29]{Ka95} for a similar result).

\begin{lemma}\label{l:operator-space}
	Let $X$, $Y$ be two Banach spaces. Let $A\in\Bb(X,Y)$ be a bounded operator with bounded $C=A^{-1}|_{AM}:AM\to M$, $B\in\Bb(X,Y)$ be a bounded operator. Let $M$, $N$ be two closed linear subspace of $X$.
	Set
	\begin{align*}
		\kappa\ :=\|C\|^{-1}(1-\delta(N,M))-\|A\|
		\delta(N,M)-\|A-B\|.
	\end{align*}	
	Then the following hold. 
	\begin{itemize}
		\item[(a)] The space $AM$ is closed, and we have
		\begin{align}\label{e:operator-space}
			\delta(AM,BN)\le \|C\|(\|A-B\|+\|B\|\delta(M,N)).
		\end{align}
		\item[(b)] Set $D_1=B^{-1}|_{BM}:BM\leadsto M$. Assume that $\|A-B\|<\|C\|^{-1}$. Then $D_1$ is an bounded operator with $\|D_1\|\le(\|C\|^{-1}-\|A-B\|)^{-1}$ and $BM$ is closed.
		\item[(c)] Assume that $\kappa>0$. Then $D=B^{-1}|_{BN}:BN\to N$ is a bounded operator with $\|D\|\le \kappa^{-1}$, $BN$ is closed and $\delta(N,M)\le (\|A\|\|C\|+1)^{-1}\le\frac{1}{2}$.
		\item[(d)] Assume that $B^{-1}$ is a bounded operator (here we do not assume that $B$ is surjective). Then we have
		\begin{align}\label{e:operator-space-reverse}
			\delta(M,N)\le\|B^{-1}\|\|(\|A-B\|+\|A\|\delta(AM,BN)).
		\end{align}
	\end{itemize}	
\end{lemma}

\begin{proof} (a)
	Let $\{x_n;\;n\in \N\}\subset M$ be a sequence of $M$ such that $\lim\limits_{n\to+\infty}Ax_n=y$. Since $C=A^{-1}|_{AM}:AM\to M$ is bounded, $\{x_n;\;n\in \N\}\subset M$ is a Cauchy sequence of $M$. Since $M$ is closed, there is an $x\in M$ such that $\lim\limits_{n\to+\infty}x_n=x$. Since $A\in\Bb(X,Y)$ is a bounded operator, we have $Ax=y$. Then $AM$ is closed.	
	
	If $M=\{0\}$, we have $\delta(M,N)=0$ and \eqref{e:operator-space} holds.
	
	Now we assume $M\ne\{0\}$. For each $x\in M\setminus\{0\}$ we have      
	\begin{align*}
		\dist(x,N)\le \delta(M,N)\|x\|.
	\end{align*}
	Then we have 
	\begin{align*}
		\dist(Ax,BN)&\le \|Ax-Bx\|+\dist(Bx,BN)\\
		&\le \|Ax-Bx\|+\|B\|\dist(x,N)\\
		&\le \|A-B\|\|x\|+\|B\|\|x\|\delta(M,N).
	\end{align*}
	Thus we obtain
	\begin{align*}
		\delta(AM,BN)&=\sup_{x\in M\setminus\{0\}}\frac{\dist(Ax,BN)}{\|Ax\|} \\   	
		&\le\|C\|(\|A-B\|+\|B\|\delta(M,N)).
	\end{align*}
	\newline (b) For each $x\in X$, we have
	\begin{align*}
		\|Bx\|\ge&\|Ax\|-\|(A-B)x\|\ge\|C\|^{-1}\|x\|-\|A-B\|\|x\|\\
		=&(\|C\|^{-1}-\|A-B\|)\|x\|.
	\end{align*} 
	Since $\|C\|^{-1}-\|A-B\|>0$,
	$D_1$ is an bounded operator and we have $\|D_1\|\le(\|C\|^{-1}-\|A-B\|)^{-1}$.
	By (a), $BM$ is closed.
	\newline (c) Since $\kappa>0$, we have $C\ne 0$ and $M\ne \{0\}$. Then we have $\|A\|\|C\|\ge 1$ and $\delta(N,M)\le (\|A\|\|C\|+1)^{-1}\le\frac{1}{2}$. 
	
	If $N=\{0\}$, we have $D=0$ and $BN=\{0\}$.
	
	Now we assume that $N\ne\{0\}$.  For each $x_2\in N\setminus\{0\}$, we have
	\begin{align*}
		\dist(x_2,M)\le\delta(N,M)\|x_2\|.
	\end{align*}
	Then for each $\varepsilon>0$, there exists an $x_1\in M$ such that 
	\begin{align*}
		\|x_2-x_1\|\le(\delta(N,M)+\varepsilon)\|x_2\|.
	\end{align*}
	Then we have 
	\begin{align*}
		(1-\delta(N,M)-\epsilon)\|x_2\|\le\|x_1\|\le(1+\delta(N,M)+\varepsilon)\|x_2\|.
	\end{align*}
	Consequently we have 
	\begin{align*}
		\|Bx_2\|&\ge\|Ax_1\|-\|Ax_2-Ax_1\|-\|Ax_2-Bx_2\|\\
		&\ge \|C\|^{-1}\|x_1\|-\|A\|\|x_2-x_1\|-\|A-B\|\|x_2\|\\ &\ge \left(\|C\|^{-1}(1-\delta(N,M)-\varepsilon)-\|A\|(\delta(N,M)+\varepsilon)-\|A-B\|\right)\|x_2\|.
	\end{align*}
	On letting $\varepsilon\to 0+$, we have $\|Bx_2\|\ge\kappa\|x_2\|$. Then $D$ is an operator and we obtain $\|D\|\le\kappa^{-1}$. By (a), the space $BN$ is closed.
	\newline (d) If $M=\{0\}$, we have $\delta(M,N)=0$ and \eqref{e:operator-space-reverse} holds.
	
	Now we assume $M\ne\{0\}$. For each $x\in M\setminus\{0\}$ we have      
	\begin{align*}
		\dist(Ax,BN)\le \delta(AM,BN)\|Ax\|.
	\end{align*}
	Then we have 
	\begin{align*}
		\dist(x,N)&\le \|B^{-1}\|\|\dist(Bx,BN)\\
		&\le \|B^{-1}\|\|(\|Ax-Bx\|+\dist(Ax,BN))\\
		&\le \|B^{-1}\|\|(\|A-B\|\|x\|+\|A\|\|x\|\delta(AM,BN)).
	\end{align*}
	Thus we obtain
	\begin{align*}
		\delta(M,N)&=\sup_{x\in M\setminus\{0\}}\frac{\dist(x,N)}{\|x\|} \\   	
		&\le\|B^{-1}\|\|(\|A-B\|+\|A\|\delta(AM,BN)).
	\end{align*}
\end{proof}

We have the following key fact. 

\begin{lemma}\label{l:strong-symplectic-by-maximal-isotropic}
	Let $(X,\omega)$ be a strong symplectic Banach space with a closed isotropic subspace $\lambda$. Then the following hold.
	\begin{itemize}
		\item [(a)] The symplectic reduction $(\lambda^\omega/\lambda, \wt{\omega})$ is a strong symplectic Banach space.
		\item[(b)] Assume that $\lambda$ is a maximal isotropic subspace with $\lambda\subsetneq\lambda^\omega$. Then we have $\gamma_\lambda>0$.
	\end{itemize}	
\end{lemma}

\begin{proof}
	(a) By \cite[Lemma 1.2.6]{BoZh18}, we have $\lambda=\lambda^{\omega\omega}$. Let $f\in(\lambda^\omega/\lambda)^*$ be a bounded semi-linear form. Then we have $f\circ\pi_{\lambda}\in (\lambda^\omega)^*$. By Hahn-Banach theorem, there is an $F\in X^*$ such that $F|_{\lambda^\omega}=f\circ\pi_{\lambda}$. Since 
	$(X,\omega)$ is a strong symplectic Banach space, there is an $x\in X$ such that $\omega(x,y)=Fy$ for all $y\in X$.
	Then we have 
	\begin{align*}
		\omega(x,y)=Fy=(f\circ\pi_{\lambda})y=0
	\end{align*}
	for all $y\in\lambda$. Consequently, we have $x\in\lambda^\omega$. Then we obtain $f=L_{\tilde\omega}(x+\lambda)$. Since $f$ is arbitrary, we obtain that $(\lambda^\omega/\lambda, \wt{\omega})$ is a strong symplectic Banach space. 
	\newline (b) By Proposition \ref{p:maximal-isotropic-properties}.c, there is an $h=\pm 1$ such that $\mi h_\lambda\tilde\omega$ is positive definite on $\lambda^\omega/\lambda$. Then \eqref{e:delta-lambda} follows from (a).
\end{proof}

\begin{theorem}\label{t:stable-max-isotropic-strong}
	Let $(X,\omega_0)$ be a strong symplectic Banach space with a maximal isotropic subspace $\lambda$. 
	Let $(X,\omega)$ be a symplectic Banach space with a closed isotropic subspace of $\mu$. 
	We set 
	\begin{align*}
		a\ :&=\|L_{\omega_0}^{-1}\|\|L_{\omega_0}\|, \\
		b(\omega)\ :&=\|L_{\omega_0}^{-1}\|\|L_\omega-L_{\omega_0}\|,\\
		f(\mu,\omega)\ :&=b(\omega)+(a+b(\omega))\delta(\lambda,\mu),\\
		g(\mu,\omega)\ :&=			    
			\left(|h_\lambda|\gamma_\lambda^{-1}\left(\|L_\omega-L_{\omega_0}\|+\|L_{\omega_0}\|f(\mu,\omega)(2+f(\mu,\omega))\right)\right)^{\frac{1}{2}}. 				
	\end{align*}
	
	Assume that there hold that 
	\begin{align}\label{e:stable-max-isotropic-strong}
		f(\mu,\omega)+g(\mu,\omega)
		<\frac{1-\delta(\lambda,\mu)}{1+\delta(\lambda,\mu)}.
	\end{align}
	Then $\mu$ is a maximal isotropic subspace of the strong symplectic Banach space $(X,\omega)$ with $h_\mu\in\{h_\lambda,0\}$.
\end{theorem}

\begin{proof} 	
	Since $L_{\omega_0}$ is an invertible skew-adjoint operator, by Proposition \ref{p:index-selfadjoint}, $X$ is reflexive. Since $b(\omega)<1$, $L_\omega$ has bounded inverse. Then $(X,\omega)$ is a strong symplectic Banach space, and we have  $\lambda^{\omega_0}=(L_{\omega_0}\lambda)^\bot$ and
	$\mu^\omega=(L_\omega\mu)^\bot$. By Lemma \ref{l:strong-symplectic-by-maximal-isotropic}, $\gamma_\lambda>0$ holds if $h_\lambda\ne 0$.
	
    By Lemma \ref{l:operator-space} and \cite[Theorem IV.2.9]{Ka95}, we have
	\begin{align*}
		\delta(\mu^\omega,\lambda^{\omega_0})\le f(\mu,\omega).
	\end{align*}  
	Then we have 
	\begin{align*}
		l(\mu,\omega)\le g(\mu,\omega).
	\end{align*}    
    By Theorem \ref{t:stable-max-isotropic} and \eqref{e:stable-max-isotropic-strong}, $\mu$ is a maximal isotropic subspace of the strong symplectic Banach space $(X,\omega)$ with $h_\mu\in\{h_\lambda,0\}$.
\end{proof}

\begin{theorem}\label{t:family-max-isotropic-strong}
	Let $X$ be a Banach space with a continuous family of strong symplectic structure $\{\omega(s)\}_{s\in J}$, where $J$ is a connected topological space. Let $\{\lambda(s)\}_{s\in J}$ be a continuous family of closed subspace of $X$ such that $\lambda(s)\subset \lambda(s)^{\omega(s)}$ holds for each $s\in J$. 
	Assume that there is an $s_0\in J$ such that $\lambda(s_0)$ is a maximal isotropic subspace of $(X,\omega(s_0))$. Then $\lambda(s)$ is a maximal isotropic subspace of $(X,\omega(s))$ with $h_{\lambda(s)}=h_{\lambda(s_0)}$ for each $s\in J$. 
\end{theorem}

\begin{proof} 
	By \cite[Lemma 1.2.6]{BoZh18}, we have $\lambda(s)= \lambda(s)^{\omega(s)\omega(s)}$ for each $s\in J$. 
	By Lemma \ref{l:strong-symplectic-by-maximal-isotropic}, we have $\gamma_{\lambda(s_0)}>0$ if $h_{\lambda(s_0)}\ne 0$, and each $s\in J$ such that $\lambda(s)$ is a maximal isotropic subspace of $(X,\omega(s))$ satisfying $\gamma_{\lambda(s)}>0$ if $h_{\lambda(s)}\ne 0$. Then our result follows from Theorem \ref{t:family-max-isotropic}.
\end{proof}

\section{Stabilities with relative bounded perturbations}\label{s:stability-relative-bounded}

Although Theorem \ref{t:stable-max-isotropic} is rather general, it is not very convenient
for application since the definition of the gap $\delta(\lambda, \mu)$ is complicated.
A less general but more convenient criterion is furnished by relatively
bounded perturbations.

In this section we will prove similar stability results in \cite[Sections V.4.1, V.4.2]{Ka95} in our setting. Firstly we introduce the notion of relative boundedness (cf. \cite[Section IV.1.1]{Ka95}, \cite[Theorem III.7.4]{Cr98}).

\begin{definition}\label{d:relative-bounded}
	Let $X$, $Y$ be two Banach spaces and $T\colon X\leadsto Y$ be a linear relation. A linear relation $A\colon X\leadsto Y$ is called to be {\em relative bounded with respect to $T$} or simply {\em $T$-bounded} if $\dom(T)\subset\dom(A)$, $\overline{T0}\supset A0$, and for some $p\ge 1$, there are non-negative constants $a_p$, $b_p$ such that
	\begin{align}\label{e:relative-bounded}
		\inf\limits_{v\in Ax}\|v\|\le \left(a_p^p\|x\|^p+b_p^p\|y\|^p\right)^\frac{1}{p}.
	\end{align}
	for each $x\in\dom(T)$ and $y\in Tx$. The greatest lower bound $b_0$
	of all possible constants $b_p$ in \eqref{e:relative-bounded} will be called the {\em relative bound of $A$
		with respect to $T$} or simply the {\em $T$ bound of $A$}. 
\end{definition}

As in \cite[Section V.4.1]{Ka95}, by H\"older's inequality, we have
\begin{align*}
	\varepsilon a_p\|x\|+(1-\varepsilon^p)^\frac{1}{p}b_p\|y\|\le \left(a_p^p\|x\|^p+b_p^p\|y\|^p\right)^\frac{1}{p}\le a_p\|x\|+b_p\|y\|
\end{align*}
for $\varepsilon\in(0,1)$.
Then the notion of relative boundedness and the number $b_0$ does not depend on the choice of $p\ge1$.

We have the following (cf. the proof of \cite[Theorem V.4.4]{Ka95}). Recall \cite[Definition II.5.4]{Cr98}:

\begin{definition}\label{d:closable}
	Let $X$, $Y$ be two Banach spaces. Let $T$ be a linear relation from $X$ to $Y$. $T$ is called {\em closable} if $\bar T$ is an extension of $T$. 
\end{definition}

\begin{lemma}\label{l:relative-bounded-closure}
	Let $X$, $Y$ be two Banach spaces and $T\colon X\leadsto Y$ be a closable linear relation. Let $A\colon X\leadsto Y$ be a $T$-bounded linear relation such that \eqref{e:relative-bounded} holds for some $p\ge 1$. Then $\bar A\hat+T0\colon X\leadsto Y$ is a $\bar T$-bounded linear relation with the same constants $a_p$ and $b_p$, that is, 
	\begin{align}\label{e:relative-bounded-closure}
		\inf\limits_{v\in \bar Ax+T0}\|v\|\le \left(a_p^p\|x\|^p+b_p^p\|y\|^p\right)^\frac{1}{p}.
	\end{align}
	for each $x\in\dom(\bar T)$ and $y\in \bar Tx$. Moreover, we have
	\begin{align}\label{e:closure-sum-T-A}
		\overline{T\hat+A}\supset\bar T\hat+\bar A.
	\end{align}
\end{lemma}

\begin{proof} 
	Since $T$ is closable, the set $T0=\bar T0$ is closed in $X$.
	Denote by $\pi\colon X\times Y\to (X\times Y)/(\{0\}\times T0)$ the natural projection along $\{0\}\times T0$. By the proof of \cite[Theorem V.4.4]{Ka95}, we have 
	\begin{align*}
		 \|\pi(\bar A)x\|&\le \left(a_p^p\|x\|^p+b_p^p\|\pi(\bar T)x\|^p\right)^\frac{1}{p},\\
		\overline{\pi(T)\hat+\pi(A)}&\supset\overline{\pi(T)}\hat+\overline{\pi(A)}.
	\end{align*} 
    Since $\overline{T0}\supset A0$, our lemma follows. 
\end{proof}

We need a stability result for the closeness. 

\begin{proposition}\label{p:stability-closeness}
	Let $X$, $Y$ be two Banach spaces. Let $T$, $S$ be linear relations from $X$ to $Y$ such that $S0=T0=W$, and
	\begin{align}\label{e:relative-bounded-T-S}
		\inf\limits_{v\in(T-S)x}\|v\|\le a\|x\|+ b^{\prime}\inf\limits_{y\in Tx}\|y\| + b^{\prime\prime}\inf\limits_{y\in Sx}\|y\|  
	\end{align} 
	for all $x\in\dom(T)=\dom(S)=D$, where $a$, $b^{\prime}$, $b$ are non-negative constants and $b^{\prime} < 1$, $b^{\prime\prime} < 1$. Set $b\ :=\max\{b^{\prime},b^{\prime\prime}\}$ and $T(\kappa)\ := T+\kappa (T-S)$ for all $\kappa\in[0,1]$. Then $S$  is
	closable if and only if $T$ is closable; in this case the closures of $T$ and $S$ have
	the same domain, and we have
	\begin{align}\label{e:delta-bar-T-S}
		\hat\delta(\overline{T(\kappa^\prime)},\overline{T(\kappa)})\le \frac{|\kappa^\prime-\kappa|(a^2+(b^{\prime}+b^{\prime\prime})^2)^\frac{1}{2}}{1-b-|\kappa^\prime-\kappa|(b^\prime+b^{\prime\prime})}
	\end{align}
	for all $\kappa^\prime, \kappa\in[0,1]$ with $|\kappa^\prime-\kappa|(b^\prime+b^{\prime\prime})<1-b$. In particular $S$ is closed if and only if $T$ is.
\end{proposition}

\begin{proof}
	1. Denote by $\pi\colon X\times Y\to (X\times Y)/(\{0\}\times W)$ the natural projection along $\{0\}\times W$. Either $S$ or $T$ is
	closable will get that $W$ is closed. Then $X$ and $Y/W$ are Banach spaces, and $\pi(T)$ and $\pi(S)$ are linear operators from $X$ to $Y/W$ with the same domain $D$. By \eqref{e:relative-bounded-T-S}, we have
	\begin{align}\label{e:pi-relative-bounded-T-S}
		\|\pi(T-S)x\|\le a\|x\|+ b^{\prime}\|\pi(T)x\| + b^{\prime\prime}\|\pi(S)x\|
	\end{align}  
	for all $x\in D$. By \cite[Theorem IV.1.3]{Ka95}, $S$  is
	closable if and only if $T$ is closable; in this case the closures of $T$ and $S$ have
	the same domain. In particular $S$ is closed if and only if $T$ is. 
	\newline 2.
	Assume that both $T$ and $S$ are closable. By \eqref{e:pi-relative-bounded-T-S}, we have
	\begin{align}\label{e:pi-bar-relative-bounded-T-S}
		\|\pi(\bar T-\bar S)x\|\le a\|x\|+ b^{\prime}\|\pi(\bar T)x\| + b^{\prime\prime}\|\pi(\bar S)x\|
	\end{align}  
	for all $x\in \dom(\bar T)=\dom(\bar S)$. 
	Set $A\ :=\pi(\bar T-\bar S)$. By \cite[(IV.1.5)]{Ka95}, we have 
	\begin{align}\label{e:pi-bar-relative-bounded-T-A}
		\|Ax\|\le \frac{1}{1-b}\left(a\|x\|+ (b^{\prime}+b^{\prime\prime})\|\pi(\overline{T(\kappa)})x\|\right)
	\end{align} 
	for all $x\in \dom(\bar T)=\dom(\bar S)$ and $\kappa\in[0,1]$.
	By \cite[Lemma A.3.1.d]{BoZh18} and \cite[Theorem IV.2.14]{Ka95}, we have
	\begin{align*}
		\hat\delta(\overline{T(\kappa^\prime)},\overline{T(\kappa)})&=\hat\delta(\pi(\overline{T(\kappa^\prime)}),\pi(\overline{T(\kappa)}))\\
		&\le \frac{|\kappa^\prime-\kappa|(a^2+(b^{\prime}+b^{\prime\prime})^2)^\frac{1}{2}}{1-b-|\kappa^\prime-\kappa|(b^\prime+b^{\prime\prime})}
	\end{align*}
	for all $\kappa^\prime, \kappa\in[0,1]$ with $|\kappa^\prime-\kappa|(b^\prime+b^{\prime\prime})<1-b$.
\end{proof}

We introduce the notions of essentially selfadjoint and essentially maximal symmetric as follows.

\begin{definition}\label{d:essentially-selfadjoint}
	Let $X$ and $Y$ be two Banach spaces over $\KK$. Let $\Omega\colon X\times Y\to\KK$ be a non-degenerate sesquilinear form. Let $\omega$ be the associated symplectic form of $\Omega$ defined by \eqref{e:natural-symplectic-structure}.
	Let $A\colon X\leadsto Y$ be a linear relation. We call $A$ {\em essentially selfadjoint} ({\em essentially maximal symmetric}, respectively) with respect to $\Omega$ if $A$ is closable and $\bar A$ is selfadjoint (maximal isotropic,respectively).
\end{definition}

Then we have the following similar result of \cite[Theorem V.4.5]{Ka95}. Note that the special case of Hilbert space situation \cite[Theorems V.4.3 and V.4.4]{Ka95} is due to F. Rellich \cite{Re39III}. 

\begin{theorem}\label{t:stable-essentially-selfadjoint-maximal-symmetric}
	Let $X$ and $Y$ be two Banach spaces over $\KK$. Let $\Omega\colon X\times Y\to\KK$ be a non-degenerate sesquilinear form with $L_\Omega^{-1}\in\Bb(Y^*,X)$. Let $T$, $S$ be two symmetric relations with respect to $\Omega$ such that $\dom(T)=\dom(S) = D$, $T0=S0$ and
	\begin{align}\label{e:relative-bounded-T-S-symmetric}
		\inf\limits_{v\in(T-S)x}\|v\|\le a\|x\|+ b(\inf\limits_{y\in Tx}\|y\| + \inf\limits_{y\in Sx}\|y\|)  
	\end{align} 
	holds for all $x\in D$, where $a$, $b$ are non-negative constants and $b< 1$.
	Then $S$ is essentially
	selfadjoint (or maximal symmetric) with respect to $\Omega$ if and only if $T$ is; in this case, $\bar S$ and $\bar T$ have the same domain.
	In particular, $S$ is selfadjoint (or maximal symmetric) with respect to $\Omega$ if and only if $T$ is.
\end{theorem}

\begin{proof} 
	Since $L_\Omega^{-1}\in\Bb(Y^*,X)$, the associated symplectic form of $\Omega$ defined by \eqref{e:natural-symplectic-structure} is strong. 
	By Proposition \ref{p:stability-closeness}, we can replace $S$ and $T$ by $\bar S$ and $\bar T$. Then our results follows from Proposition \ref{p:stability-closeness} and Theorem \ref{t:family-max-isotropic-strong}.
\end{proof}

\section{Skew-adjoint relations between real Banach spaces}\label{s:real-skew-adjoint}

\subsection{Skew-adjoint relations between real Hilbert spaces}\label{ss:real-skew-adjoint-Hilbert}

We denote by $V_{\C}\ :=V\otimes_{\R}\C$ for each real vector space $V$. 
Let $(X,\lla\cdot,\cdot\rra)$ be a real Hilbert space. Then $X^2$ is a Hilbert space with the inner product defined by
\begin{align*}
	\lla(x_1,y_1),(x_2,y_2)\rra=\lla x_1,x_2\rra+\lla y_1,y_2\rra
\end{align*}
for each $x_1,y_1,x_2,y_2\in X$.
Set $Q(x,y)\;:=\lla x,y\rra$ for each $x,y\in X$. Denote by $I_X$ the identity map on $X$, and 
\begin{align*}
	J_{X^2}\: =\left(\begin{array}{cc}
		0	&-I_X  \\
		I_X	& 0
	\end{array}\right).
\end{align*} 
Then we have $\omega(u,v)=\lla J_{X^2}u,v\rra$ for all $u,v\in X^2$, where $\omega$ is the symplectic structure of $X^2$ associated with $Q$. We denote by $\U(X_\C)$ the unitary group on $X_\C$ and $\OO(X)$ the orthogonal group on $X$. We have the following structure of $\Rr^{-1,sa}(X)$.

\begin{proposition}\label{p:real-skew-adjoint-Hilbert}
	Let $(X,Q=\lla\cdot,\cdot\rra)$ be a real Hilbert space and $Y$ be a real linear space. Let $\Omega\colon X\times Y$ be a bilinear form such that $R_\Omega\colon Y\to X^*$ is a real linear isomorphism. Then we have a homemorphism
	$\varphi_\Omega\colon\U(X_\C)\to\Rr^{-1,sa}(\Omega_\C)$ defined by
	\begin{align}\label{e:skew-adjoint-Hilbert}
		\varphi_\Omega(U)=R_\Omega^{-1}\circ R_Q\circ (I_{X_C}-U)\circ(I_{X_C}+U)^{-1}
	\end{align}
    for all $U\in \U(X_C)$ such that $\varphi_\Omega(\OO(X))=\{A_\C;\;A\in\Rr^{-1,sa}(\Omega)\}$. For each $T=\varphi_\Omega(U)\in\Rr^{-1,sa}(\Omega_\C)$, we have $\ker_\C T=\ker_\C(U-I)$ and $T_\C0=(R_\Omega^{-1}\circ R_Q)\ker_\C(U+I)$.
\end{proposition} 

\begin{proof} 
	1. Let $A\colon X\leadsto X$ be a linear relation. Then we have
	\begin{align*}
		(R_\Omega^{-1}\circ R_Q\circ A)^\Omega=R_\Omega^{-1}\circ (R_Q\circ A)^*=R_\Omega^{-1}\circ R_Q\circ A^Q.
	\end{align*} 
	\newline 2. Denote by $X^{\pm}\ :=\{(x,\pm\mi x);\;x\in X_\C\}$, Then we have $J_{X^2}u=\mp \mi u$ for all $u\in X^{\pm}$. By \cite[Proposition 2]{BoZh13} we have 
	\begin{align*}
		\Rr^{1,sa}(X_\C)&=\{\{(x,\mi x)+(Ux, -\mi Ux);\;x\in X_\C\};\;U\in\U(X_\C)\}\\
		&=\{\mi (I_{X_C}-U)\circ(I_{X_C}+U)^{-1};\;U\in\U(X_\C)\}.
	\end{align*}
    Then we have
    \begin{align*}
    	\Rr^{-1,sa}(X_\C)&=\{\{(x, x)+(Ux, - Ux);\;x\in X_\C\};\;U\in\U(X_\C)\}\\
    	&=\{ (I_{X_C}-U)\circ(I_{X_C}+U)^{-1};\;U\in\U(X_\C)\}.
    \end{align*}    
    By Step 1, $\varphi_{\Omega_X}$ defined by \eqref{e:skew-adjoint-Hilbert} is a homemorphism and $\varphi_{\Omega_X}(\OO(X))\subset\{A_\C;\;A\in\Rr^{-1,sa}(X)\}$.
    \newline 3. Let $A\in\{A_\C;\;A\in\Rr^{-1,sa}(X)\}$. By Step 1, there is a $U=C+\mi D\in\U(X_\C)$, $C,D\in\Bb(X)$ such that 
    \begin{align*}
    	A_\C&=\{(x, x)+(Ux, - Ux);\;x\in X_\C\}\\
    	&=\{(u(x,y),v(x,y));\;x,y\in X\}.
    \end{align*}
    where $u(x,y)=x+Cx-Dy+\mi(y+Cy+Dx)$ and $v(x,y)=x-Cx+Dy+\mi(y-Cy-Dx)$. From this we obtain
    \begin{align*}
    	A&=\{(x+Cx-Dy,x-Cx+Dy);\;x,y\in X\}\\
    	&=\{(y+Cy+Dx,y-Cy-Dx);\;x,y\in X\}\\
    	&\subset A_\C.
    \end{align*}
    We conclude that $Cx-Dy=Ux$ for all $x,y\in X$. Then we have $D=0$ and $U=C\in\OO(X)$.
    \newline 4. By Step 1, Step 2 and Step 3, the map $\varphi_\Omega$ defined by \eqref{e:skew-adjoint-Hilbert} satisfies the desired properties.
\end{proof}

\subsection{Stability theorem of mod 2 index}\label{ss:relative-dimension}

In this subsection we prove Theorem \ref{t:stable-mod-2-index}. 

\begin{proof}[Proof of Theorem \ref{t:stable-mod-2-index}] 
	We divide the proof into five steps. We use the same symbol to denote the complexified sesquilinear forms of the corresponding real bilinear form.
	\newline 1. Since $T\colon X\leadsto Y$ is a skew-adjoint linear Fredholm relation with respect to $\Omega_0$ with index $0$, there is a real linear subspace $V$ of $Y$ such that 
	\begin{align}
		\label{e:image-split-T}		
		Y=\image T\oplus V,\quad
		\dim\ker T=\dim V<+\infty.
	\end{align}      
    \newline 2. Assume that $\dim X<+\infty$. Define $Q_T\colon\dom (T)\times\dom(T)\to\R$ by $Q_T(x,y)\ :=\Omega_0(x,Ty)$ for all $x,y\in\dom(T)$. Since $T$ is skew-adjoint, we have $\dim T=\dim X=\dim Y$. By Witt decomposition theorem, we have 
    \begin{align}\label{e:mod-2-index-finite-dimension}
    	\dim\dom(T)=\dim T-\dim T0=2m^-(\mi Q_T)+\dim\ker Q_T. 
    \end{align}
    By Proposition \ref{p:maximal-isotropic-properties}, we have $\dim\ker Q_T=\dim\ker T$. Then we have
    \begin{align}\label{e:mod-2-index-finite-dimension-1}
    	\dim\dom(T)=\dim X-\dim T0=2m^-(\mi Q_T)+\dim\ker T. 
    \end{align}
    \newline 3. Set $Z\ :=X_{\C}\times Y_{\C}$. 
    We denote by $\omega_0$ and $\omega$ the associated symplectic structure of $\Omega_0$ and $\Omega$ on $Z$ respectively. We set
    \begin{align*}
    	\alpha\ :=X_{\C}\times\{0\},\quad
    	\beta\ :=\{0\}\times V_{\C},\quad
    	\gamma\ :=(\mi\circ S)_{\C}, \quad
    	\lambda\ :=(\mi T)_{\C}.
    \end{align*}
    Then $\alpha$ is a Lagrangian subspace of $(Z,\omega_0)$ and $(Z,\omega)$, $\beta$ is a finite dimensional isotropic subspace of $(Z,\omega_0)$ and $(Z,\omega)$, $\gamma$ is an isotropic subspace of $(Z,\omega)$, and $\lambda$ is a Lagrangian subspace of $(Z,\omega_0)$. 
    By \cite[Proposition 1.3.3]{BoZh18} we have 
    \begin{align*}
    	&Z=\alpha\cap\beta^{\omega_0}\oplus\beta\oplus\lambda=\alpha+\beta^{\omega_0}.
    \end{align*}
    We have 
    $(\alpha\cap\beta^\omega)^\omega\supset\alpha+\beta$.
    By \cite[Lemma 1.1.2]{BoZh18}, we have
    \begin{align*}
    	\dim_{\C}(\alpha\cap\beta^\omega)^\omega/\alpha\le&\dim_\C\alpha/(\alpha\cap\beta^\omega)\\
    	=&\dim_\C X_\C/V_\C^{\Omega,l}\\
    	=&\dim_\C V_\C=\dim_\C \beta.
    \end{align*} 
    Then we have
    \begin{align*}
    	(\alpha\cap\beta^\omega)^\omega=\alpha+\beta.
    \end{align*} 
    By \cite[Proposition A.3.13 and Lemma 3.1.1]{BoZh18}, there is a $\delta>0$ such that for each closed skew symmetric linear relation $S$ with respect to $\Omega$ with $\hat\delta(S,T)+\|L_\Omega-L_{\Omega_0}\|<\delta$, we have
    \begin{align*}
    	\{0\}&=\beta\cap\gamma,\\
    	Z&=\alpha\cap\beta^\omega+\beta+\gamma=\alpha+\beta^\omega,\\
    	\{0\}&=Z^\omega=(\alpha+\beta+\gamma)^\omega=\alpha\cap\beta^\omega\cap\gamma,
    \end{align*}   
    By \cite[Theorem IV.4.30 and Problem IV.4.6]{Ka95}, we have 
    \begin{align*}
    	0=\Index T=\Index S= \Index_\C(\alpha+\beta,\gamma)-\dim_\C \beta.
    \end{align*}   
	By Lemma \ref{l:selfadjoint-index-non-positive}, $S$ is skew-adjoint with respect to $\Omega$. 
	\newline 4.  
    We define the sesquilinear form $Q$ on $\alpha\cap(\beta+\gamma)$ by $Q(x_1,x_2)=\omega(x_1,y_2)$ if $x_j\in \alpha$, $y_j\in \beta$, $x_j+y_j\in\gamma$, $j=1,2$. By \cite[\S3.1]{ZWZ18},
    the form $Q$ is a well-defined symmetric form.
        
    Let $\alpha_1$ be a closed complex linear subspace of $\alpha$ with $\dim_\C\alpha/\alpha_1<+\infty$. By \cite[Problem IV.4.6]{Ka95} we have
    \begin{align*}
    	\dim_\C(\gamma\cap(\alpha_1+\beta))
    	&=\Index_\C(\alpha_1+\beta,\gamma)+\dim_\C Z/(\alpha_1+\beta+\gamma)\\
    	&=\Index_\C(\alpha_1,\gamma)+\dim_\C\beta+\dim_\C Z/(\alpha_1+\beta+\gamma)\\
    	&=\dim_\C\beta+\dim_\C Z/(\alpha_1+\beta+\gamma)-\dim_\C\alpha/\alpha_1.
    \end{align*}

     Since $\alpha_1\cap\beta=\beta\cap\gamma=\{0\}$, there is a linear isomorphism $\phi\colon\gamma\cap(\alpha_1+\beta)\to\alpha_1\cap(\beta+\gamma)$ defined by $z-\phi(z)\in\beta$ for each $z\in\gamma\cap(\alpha_1+\beta)$. Then we have 
     \begin{align}\label{e:dimension-domain}
    	 \begin{split}
    		 \dim_\C(\alpha_1\cap(\beta+\gamma))=\dim_\C\beta+\dim_\C Z/(\alpha_1+\beta+\gamma)-\dim_\C\alpha/\alpha_1.
    	\end{split}    	
    \end{align}
    Since $\alpha+\beta+\gamma=Z$, we have
    \begin{align}\label{e:dimension-domain-ker}
    	\begin{split}
    		\dim_\C(\alpha\cap(\beta+\gamma))=\dim_\C\beta=\dim\ker T.
    	\end{split}    	
    \end{align}
        
    We have  $(\alpha+\beta)^\omega=\alpha\cap\beta^\omega$. Then we have 
    \begin{align*}
    	(\alpha+\beta)^{\omega\omega}=(\alpha\cap\beta^\omega)^\omega=\alpha+\beta.
    \end{align*} 
    Since $Z=\alpha+\beta^\omega$, by \cite[Lemma 1.3.2]{BoZh18}, we have 
    \begin{align*}
    	\dim_\C\alpha/(\alpha\cap\beta^\omega)=\dim_\C Z/\beta^\omega=\dim_\C\beta.
    \end{align*}
    By \cite[Problem IV.4.6]{Ka95} we have 
    \begin{align*}
    	\Index_\C(\alpha+&\beta,\gamma)+\Index_\C((\alpha+\beta)^\omega,\gamma)\\
    	=&\Index_\C(\alpha,\gamma)+\dim_\C(\alpha+\beta))/\alpha+\\
    	&\Index_\C(\alpha,\gamma)-\dim_\C\alpha/(\alpha+\beta)^\omega\\
    	=&2\Index_\C(\alpha,\gamma)=0.
    \end{align*}
    By \cite[Corollary 1.2.12]{BoZh18} we have
    \begin{align*}
    	(\gamma\cap(\alpha+\beta))^{\omega}=\gamma+(\alpha+\beta)^\omega=\gamma+\alpha\cap\beta^\omega.
    \end{align*}
    By \cite[Lemma 3.3]{ZWZ18} and \cite[Lamma A.1.1]{BoZh18} we have
    \begin{align*}
    	\ker _\C Q&=\alpha\cap(\beta\cap(\gamma\cap(\alpha+\beta))^{\omega}+\gamma)\\
    	&=\alpha\cap(\beta\cap(\gamma+\alpha\cap\beta^\omega)+\gamma)\\
    	&=\alpha\cap(\beta+\gamma)\cap(\gamma+\alpha\cap\beta^\omega)\\
    	&=(\beta+\gamma)\cap(\alpha\cap\gamma+\alpha\cap\beta^\omega).
    \end{align*}
    Since $\alpha\cap\beta^\omega+\beta+\gamma=Z$ and $\alpha\cap\beta^\omega\cap\gamma=\{0\}$,
    by \eqref{e:dimension-domain} and \cite[Problem IV.4.6]{Ka95} we have
    \begin{align*}
    	\dim_\C\ker_\C Q
    	=&\dim_\C\beta-\dim_\C\alpha/(\alpha\cap\gamma+\alpha\cap\beta^\omega)\\
    	=&\dim_\C\beta-\dim_\C\alpha/(\alpha\cap\beta^\omega)+\\
    	&\dim_\C(\alpha\cap\gamma+\alpha\cap\beta^\omega)/(\alpha\cap\beta^\omega)\\
    	=&\dim_\C(\alpha\cap\gamma)/(\alpha\cap\beta^\omega\cap\gamma)\\
    	=&\dim_\C(\alpha\cap\gamma)
    	=\dim\ker S.
    \end{align*}
    \newline 5.  
    We have $\alpha\cap(\beta+\gamma)=S_\C^{-1}V_{\C}\times\{0\}$. For all $x_1,x_2\in S_\C^{-1}V_{\C}$, we take a $y_1\in V_{\C}\cap S_\C x_1$ and a $y_2\in V_{\C}\cap S_\C x_2$. 
    Then we have 
    \begin{align*}
    	Q((x_1,0),(x_2,0))&=\omega((x_1,0), (0,\mi y_2))\\
    	&=\Omega(x_1,\mi y_2).
    \end{align*}
    Note that for $x_1,x_2\in S^{-1}V$, we can take a $y_1\in V\cap S x_1$ and a $y_2\in V\cap S x_2$.
    Then $\mi Q$ is a real skew symmetric form on the real vector space $S^{-1}V\times\{0\}$. By Step 1, Step 4,  \eqref{e:mod-2-index-finite-dimension} and \eqref{e:dimension-domain-ker} we have 
    \begin{align*}
    	\dim\ker T=&\dim_\C(\alpha\cap(\beta+\gamma))=\dim (S^{-1}V)=\dim\ker \mi Q|_{S^{-1}V\times\{0\}}\\
    	=&\dim_\C\ker_{\C} Q=\dim\ker S\mod 2.
    \end{align*}
\end{proof}

\begin{proof}[Proof of Theorem \ref{t:family-mod-2-index}]
	Set $J_a\ :=\{s\in J;\;\dim\ker T(s)=a\mod 2\}$. Then we have $J=J_0\cup J_1$. By Theorem \ref{t:stable-mod-2-index}, the sets $J_0$ and $J_1$ are open subsets of $J$. Since $J$ is connected, we have $J=J_0$ oe $J=J_1$. Then $\dim\ker T(s)\mod 2$ is a constant in $\Z_2$ for $s\in J$.
\end{proof}

\begin{proof}[Proof of Corollary \ref{c:real-symplectic-extension-even}] 
	Note that $(X\oplus\R^n)^*=X^*\oplus\R^n$. 
	Denote by $i_0:\R^n\to X\oplus\R^n$, $i_1:X\to X\oplus\R^n$, $j_0:\R^n\to X^*\oplus\R^n$, $j_1:X^*\to X^*\oplus\R^n$ the natural embeddings, and $p_0:X\oplus\R^n\to\R^n$, $p_1:X\oplus\R^n\to X$, $q_0:X^*\oplus\R^n\to\R^n$, $q_1:X^*\oplus\R^n\to X^*$ the natural projections. Then we have $L_\Omega=j_1\circ L_\Omega$ and $p_1\circ i_1=I_{X_1}$, and
	\begin{align*}
		L_{\tilde\omega}|_X-L_\omega&=L_{\tilde\omega}\circ i_1-j_1\circ L_\omega\\
		&=(L_{\tilde\Omega}-j_1\circ L_\omega\circ p_1)\circ i_1\\
		&=(L_{\tilde\omega}-j_1\circ L_\omega\circ p_1)\circ(I_{X\oplus\R^n}-i_0).
	\end{align*}
	Since $i_0$ is an operator of finite rank and the operator $L_{\tilde\omega}|_X-L_\omega$ is compact, the operator $L_{\tilde\omega}-j_1\circ L_\omega\circ p_1$ is compact. Since $L_{\tilde\omega}$ and $L_\omega$ are injective, $(X\oplus\R^n,\tilde\omega)$ is a strong symplectic Banach space, if and only if $L_{\tilde\omega}$ is an isomorphism if and only if $j_1\circ L_\omega\circ p_1$ is Fredholm of index $0$, if and only if $L_\omega$ is Fredholm of index $0$, if and only if $L_\omega$ is 
	an isomorphism, if and only if $(X,\omega)$ is a strong symplectic Banach space.
	
	In the affirmative case, we have a continuous family of skew-adjoint Fredholm operators $\{L_{\tilde\omega}+s(j_1\circ L_\omega\circ p_1-L_{\tilde\omega})\}_{s\in[0,1]}$ with index $0$. By Theorem \ref{t:family-mod-2-index}, we have
	\begin{align*}
		n=\dim\ker j_1\circ L_\omega\circ p_1
		=\dim\ker L_{\tilde\omega}
		=0\mod 2,
	\end{align*}
	i.e. $n$ is even.
\end{proof}

\subsection{Path components}\label{ss:path-components-skew}
In this subsection, we study the path components of the set of skew-adjoint Fredholm relations between real Banach spaces with index $0$ and prove Theorem \ref{t:path-components-real-skew}.

Let $X$ and $Y$ be two Banach space, and $\Omega\colon X\times Y\to\KK$ be a bounded non-degenerate sesquilinear form. Let $h$ be in $\KK$ with $|h|=1$.
For $k\in\N$, we define 
\begin{align}\label{e:skew-adjoint-k}
	\Ff\Rr_{k,0}^{h,sa}(\Omega)\ :&=\{T\in \Ff\Rr_0^{h,sa}(\Omega);\;\dim\ker T=k\}.
\end{align}

We need the following lemmas.

\begin{lemma}\label{l:transversal-real-skew}
	Let $X$ and $Y$ be two vector space over $\KK$, and $\Omega\colon X\times Y$ be a non-degenerate sesquilinear form. Let $h$ be in $\KK$ with $|h|=1$.
	Then $\Ff\Rr_{0,0}^{h,sa}(\Omega)$ is an affine space over $\KK$. 
\end{lemma}

\begin{proof}
	Since $\Ff\Rr_{0,0}^{h,sa}(\Omega)=\{T\in \Rr^{h,sa}(\Omega);\; T^{-1}\in\Hom(Y,X)\}$, the set $\Ff\Rr_{0,0}^{h,sa}(\Omega)$ is an affine space over $\KK$. 
\end{proof}

\begin{lemma}\label{l:reduce-to-transversal-case}
	Let $X$ and $Y$ be two vector space over $\KK$, and $\Omega\colon X\times Y$ be a non-degenerate sesquilinear form. Let $h$ be in $\KK$ with $|h|=1$. Then the following hold.
	\begin{itemize}		
	    \item [(a)] Let $X_0$, $X_1$ be linear subspaces of $X$ and $Y_0$, $Y_1$ be linear subspaces of $Y$ such that 
	    \begin{align}\label{e:Omega-splitting}
	    	\begin{split}
	    		&X=X_0+ X_1,\quad Y=Y_0+ Y_1,\\
	    		&\Omega(x,y)=0\quad\forall(x,y)\in
	    		(X_0\times Y_1)\cup(X_1\times Y_0).  
	    	\end{split}			
	    \end{align}
	    Then we have
	    \begin{align}\label{e:splitting-XY-01}
	    	&X=X_0\oplus X_1,\quad Y=Y_0\oplus Y_1,\\
	    	\label{e:XY-01}
	    	&X_0=Y_1^{\Omega,l},\quad X_1=Y_0^{\Omega,l},\quad Y_0=X_1^{\Omega,r},\quad Y_1=X_0^{\Omega,r},\\
	    	\label{e:skew-relation-splitting-XY-01}
	    	&\Ff\Rr^{h,sa}(\Omega)\supset\{T_0\oplus T_1;\;T_0\in\Ff\Rr^{h,sa}(\Omega|_{X_0\times Y_0}),
	    	T_1\in\Ff\Rr_{0,0}^{h,sa}(\Omega|_{X_1\times Y_1})\}.
	    \end{align}
	    \item [(b)] Let $X_0$, $Y_0$ be finite dimensional linear subspaces of $X$ and $Y$ respectively such that 
	    $Y=Y_0\oplus X_0^{\Omega,r}$.
	    Then we have $\dim X_0=\dim Y_0$ and
	    $X=X_0\oplus Y_0^{\Omega,l}$.
	    \item [(c)] Let $T\colon X\leadsto Y$ be a $h$-selfadjoint linear Fredholm relation with index $0$. We denote by $X_0\ :=\ker T$ and $Y_1\ :=\ran T$. Let $Y_0$ be a linear subspace of $Y$ such that $Y=Y_0\oplus Y_1$. We denote by $X_1\ :=Y_0^{\Omega,l}$. Denote by $Q_T$ the associated form of $T$ defined by Definition \ref{d:selfadjoint}.c. Then \eqref{e:Omega-splitting} holds, and we have 
	    \begin{align}\label{splitting-T-ker-Q}
	    	T&=T_0\oplus T_1,\quad Q_T=Q_{T_0}\oplus Q_{T_1},\quad\ker T=\ker Q_T,
	    \end{align}
	    where
	    \begin{align}
	    	\label{splitting-T0}T_0\ :&=X_0\times\{0\}\in\Ff\Rr^{h,sa}(\Omega|_{X_0\times Y_0}),\\
	    	\label{splitting-T1}T_1\ :&=T\cap(X_1\times Y_1)\in\Ff\Rr_{0,0}^{h,sa}(\Omega|_{X_1\times Y_1}).	    	
	    \end{align}
	\end{itemize}
\end{lemma}
    
\begin{proof}
	(a) Since $Y=Y_0+Y_1$ and $X_0\subset Y_1^{\Omega,l}$ hold, we have
	\begin{align*}
		X_0\cap Y_0^{\Omega,l}=X_0\cap Y_1^{\Omega,l}\cap Y_0^{\Omega,l}
		=X_0\cap(Y_0+Y_1)^{\Omega,l}=\{0\}. 
	\end{align*}
    Similarly we have $Y_0\cap X_0^{\Omega,r}=\{0\}$. Then $\Omega|_{X_0\times Y_0}$ is nondegenerate.
    Similarly $\Omega|_{X_1\times Y_1}$ is nondegenerate.
    	
	Since $Y=Y_0+Y_1$, $X_0\subset Y_1^{\Omega,l}$ and and $X_1\subset Y_0^{\Omega,l}$ hold, we have 
	\begin{align*}
		X_0\cap X_1\subset Y_0^{\Omega,l}\cap Y_1^{\Omega,l}=(Y_0+Y_1)^{\Omega,l}=\{0\}.
	\end{align*}
	Then we have $X=X_0\oplus X_1$. Similarly, we have $Y=Y_0\oplus Y_1$.
	
	Since $Y_1^{\Omega,l}\supset X_0$ and $Y_1^{\Omega,l}\cap X_1=\{0\}$, by \cite[Lemma A.1.1]{BoZh18} we have 
	\begin{align*}
		Y_1^{\Omega,l}=Y_1^{\Omega,l}\cap(X_0+X_1)=X_0+Y_1^{\Omega,l}\cap X_1=X_0.
	\end{align*}
    Similarly we get the other three identities of \eqref{e:XY-01}.
    
    The equation \eqref{e:skew-relation-splitting-XY-01} is clear.
    \newline (b) Since $\Omega$ is non-degenerate and $Y=Y_0\oplus X_0^{\Omega,r}$, by \cite[Lemma 1.1.2.c]{BoZh18}, we have
    \begin{align*}
    	\dim X/Y_0^{\Omega,l}=\dim Y_0=\dim Y/X_0^{\Omega,r}=\dim X_0.
    \end{align*}
    By \cite[Lemma 1.1.2.c]{BoZh18}, we have 
    \begin{align*}
    	X_0\cap Y_0^{\Omega,l}=(X_0^{\Omega,r}+Y_0)^{\Omega,l}=Y^{\Omega,l}=\{0\}.
    \end{align*}
    Then we conclude that $X=X_0\oplus Y_0^{\Omega,l}$.
	\newline (c) 
	By \cite[Lemma 1.1.2.c]{BoZh18} we have $\dim X_0=\dim Y_0=\dim X/X_1$. By Proposition \ref{p:ker-dom-max-symmetric}, we have $X_0=Y_1^{\Omega,l}$. Since $Y=Y_0\oplus Y_1$, we have $X_0\cap X_1=(Y_0+Y_1)^{\Omega,l}=\{0\}$. Then we have $X=X_0\oplus X_1$. Then \eqref{e:Omega-splitting} holds. By (a), we have $T=T_0\oplus T_1$, where $T_0$ and $T_1$ are defined by \eqref{splitting-T0} and \eqref{splitting-T1} respectively.
	Then we have $Q_T=Q_{T_0}\oplus Q_{T_1}$. Since $\ran T_1=Y_1$, we have $\ker Q_{T_1}=\{0\}$. Then we have 
	\begin{align*}
		\ker Q_T=\ker Q_{T_0}\oplus\ker Q_{T_1}=X_0=\ker T.
	\end{align*}
\end{proof}

We define the map $\pi\colon \Ff\Rr_{k,0}^{h,sa}(\Omega)\to G(k,X)$ by 
\begin{align}\label{e:skew-adjoint-k-fibration}
	\pi(T)=\ker T.
\end{align}
We shall prove that $(\Ff\Rr_{k,0}^{h,sa}(\Omega),\pi)$ is a fiber bundle for $|h|=1$.

\begin{lemma}\label{l:continuity-AB}
	Let $X$ and $Y$ be two Banach space over $\KK$, and $\Omega\colon X\times Y\to\KK$ be a bounded non-degenerate sesquilinear form. Let $X_0$, $Y_0$ be finite dimensional linear subspaces of $X$ and $Y$ respectively such that 
	$Y=Y_0\oplus X_0^{\Omega,r}$. For each $A\in\Bb(X_0,Y_0^{\Omega,l})$, we define $f_\Omega(A)\in\Bb(X_0^{\Omega,r},Y_0)$ by 
	\begin{align}\label{e:finite-dimensional-dual}
		\Omega(x_0,f_\Omega(A)y_1)=\Omega(Ax_0,y_1)
	\end{align}
	for each $(x_0,y_1)\in X_0\times X_0^{\Omega,r}$. Then the map $f_\Omega$ is a bounded real linear map. 
\end{lemma}

\begin{proof}
	By Lemma \ref{l:reduce-to-transversal-case}.b, we have $\dim X_0=\dim Y_0$ and $X=X_0\oplus Y_0^{\Omega,l}$. By Lemma \ref{l:reduce-to-transversal-case}.a, the sesquilinear forms $\Omega_0\ :=\Omega|_{X_0\times Y_0}$ and $\Omega_1\ :=\Omega|_{Y_0^{\Omega,l}\times X_0^{\Omega,r}}$ are non-degenerate. Since $\dim X_0=\dim Y_0<+\infty$, there is a unique $f_\Omega(A)y_1\in Y_0$ for each $x_0\in X_0$ such that \eqref{e:finite-dimensional-dual} holds for all $x_0\in X_0$. 
	
	We see from the definition that $f_\Omega(A)$ is a linear map and $f_\Omega$ is a real linear map. By \eqref{e:form-operator}, there is a constant $C>0$ such that 
	\begin{align*}
		|R_{\Omega_0}f_\Omega(A)y_1(x_0)|=|\Omega(x_0,f_\Omega(A)y_1)|=|\Omega(Ax_0,y_1)|\le C\|A\|\|x_0\|\|y_1\|				
	\end{align*}
	for each $(x_0,y_1)\in X_0\times X_0^{\Omega,r}$.
	Since $\dim X_0=\dim Y_0<+\infty$, $(R_{\Omega_0})^{-1}\in\Bb(Y_0^*,X_0)$. is an bounded operator. Then we have
	\begin{align*}
		\|f_\Omega(A)\|=&\sup_{y_1\in X_0^{\Omega,r},\|y_1\|=1}\|f_\Omega(A)y_1\|\\
		\le &\sup_{y_1\in X_0^{\Omega,r}, \|y_1\|=1}\|(R_{\Omega_0})^{-1}\||R_{\Omega_0}f_\Omega(A)y_1|\\
		\le &\sup_{y_1\in X_0^{\Omega,r}, \|y_1\|=1}\|(R_{\Omega_0})^{-1}\|C\|A\|\|y_1\|\\
		=&C\|(R_{\Omega_0})^{-1}\|\|A\|.
	\end{align*}
	Thus the map $f_\Omega$ is a bounded real linear map. 	
\end{proof}

\begin{rem}\label{r:f-Omega-not-surjective}
	The map $f_\Omega$ defined in Lemma \ref{l:continuity-AB} is not surjective in general.
\end{rem}

\begin{proposition}\label{l:skew-adjoint-k-fibration}
	Let $X$ and $Y$ be two Banach space over $\KK$, and $\Omega\colon X\times Y\to\KK$ be a bounded non-degenerate sesquilinear form. 	
	Denote by $G(k,X)$ the Grassmannian of the linear subspace of $X$ of dimension $k$.
	Let $\pi$ be the map defined by \eqref{e:skew-adjoint-k-fibration}.
	Then $(\Ff\Rr_{k,0}^{h,sa}(\Omega),\pi)$ is a fiber bundle with fiber $\pi^{-1}(X_0)\simeq\Ff\Rr_{0,0}^{h,sa}(\Omega|_{X_1\times Y_1})$ at $X_0$, where $Y_1\ :=X_0^{\Omega,r}$, $Y_0$ is a finite dimensional linear subspace of $Y$ with $Y=Y_0\oplus Y_1$, and $X_1\ :=Y_0^{\Omega,l}$.  
\end{proposition}

\begin{proof} By \cite[Corollary A.3.6]{BoZh18}, $\pi$ is a continuous map.
	
	By Lemma \ref{l:reduce-to-transversal-case}. b, we have $\dim Y_0=\dim X_0=k$ and $X=X_0\oplus X_1$. By Lemma \ref{l:reduce-to-transversal-case}, $T\in \pi^{-1}(X_0)$ holds if and only if $T=T_0\oplus T_1$, where $T_0\ :=X_0\times\{0\}$ and $T_1\ :=T\cap(X_1\times Y_1)\in\Ff\Rr_{0,0}^{h,sa}(\Omega|_{X_1\times Y_1})$. Then we have    	 $\pi^{-1}(X_0)\simeq\Rr_{0,0}^{h,sa}(\Omega|_{X_1\times Y_1})$.
		
	Set $G(X,X_1)\ :=\{M\in\Ss(X);X=M\oplus X_1\}$
	. Then we have $X_0\in G(X,X_1)$. By \cite[(A.10)]{BoZh18}, the set $G(X,X_1)$ is an open affine space and we have
	\begin{equation}\label{e:set-H}
		G(X,X_1)=\{\Graph(A);A\in\Bb(X_0,X_1)\}.
	\end{equation}.
    
	By \cite[(A.10)]{BoZh18}, a linear subspace $M\in G(X,X_1)$ can be written as $M=\Graph(A)$ with $A\in\Bb(X_0,X_1)$.     
    By Lemma \ref{l:reduce-to-transversal-case}, we have 
	\begin{align*}
		Y_0\cap M^{\Omega,r}= X_1^{\Omega,r}\cap M^{\Omega,r}=(X_1+M)^{\Omega,r}=\{0\}.
	\end{align*}
    By \cite[Lemma 1.1.2.c]{BoZh18}, we have 
    \begin{align*}
    	\dim Y/M^{\Omega,r}=\dim M=\dim X_0=\dim Y_0=k.
    \end{align*}
    Then we have $Y=Y_0\oplus M^{\Omega,r}$.    
    
    By \cite[(A.10)]{BoZh18}, we have 
    $M^{\Omega,r}=\Graph(B)$ for a uniquely determined $B\in\Bb(Y_1,Y_0)$. Then we have 
    \begin{align}\label{e:skew-A-B}
    	\Omega(x_0,By_1)+\Omega(Ax_0,y_1)=\Omega(x_0+Ax_0,By_1+y_1)=0
    \end{align}
    for each $(x_0,y_1)\in X_0\times Y_1$.
    
    We have
    \begin{align*}
    	S_1\in&\Ff\Rr_{0,0}^{h,sa}(\Omega|_{X_1\times M^{\Omega,r}})\\
    	\Leftrightarrow &S_1=
    	(B+I_{Y_1})\circ T_1, T_1\subset X_1\times Y_1, S_1\oplus(X_1\times\{0\})=X_1\times M^{\Omega,r},\\
    	&h\circ (B+I_{Y_1})\circ T_1=((B+I_{Y_1})\circ T_1)^\Omega\\
    	\Leftrightarrow &S_1=
    	(B+I_{Y_1})\circ T_1, T_1\subset X_1\times Y_1, T_1\oplus(X_1\times\{0\})=X_1\times Y_1,\\
    	&h\circ T_1=T_1^\Omega\\
    	\Leftrightarrow &S_1=
    	(B+I_{Y_1})\circ T_1, T_1\in\Ff\Rr_{0,0}^{h,sa}(\Omega|_{X_1\times Y_1}).
    \end{align*}

    With the notations above, we define $\varphi_{X_1,Y_1}:\pi^{-1}(G(X,X_1))\to G(X,X_1)\times \Ff\Rr_{0,0}^{h,sa}(\Omega|_{X_1\times Y_1})$
    by 
    \begin{align}\label{e:skew-local-product}
    	\varphi_{X_1,Y_1}((M\times\{0\})\oplus((B+I_{Y_1})\circ T_1))=(M,T_1).
    \end{align}
    By \cite[Lemma 1.4]{Neu68}, we have $(M\times\{0\})\oplus((B+I_{Y_1})\circ T_1)$ varying continuously if and only if $M$ and $(B+I_{Y_1})\circ T_1$ varying continuously. By \eqref {e:skew-A-B}, we have $B=-f_\Omega(A)$, where $f_\Omega$ is defined by Lemma \ref{l:continuity-AB}. By Lemma \ref{l:continuity-AB}, the map $A\longmapsto B$ is continuous. By Lemma \ref{l:operator-space}, when $B$ is continuously varying, we have $(B+I_{Y_1})\circ T_1)$ varying continuously if and only if $T_1$ varying continuously. Then we have $(M\times\{0\})\oplus((B+I_{Y_1})\circ T_1)$ varying continuously if and only if $M$ and $T_1$ varying continuously.
    Thus $\varphi_{X_1,Y_1}$ is a homemorphism. By definition, $(\Ff\Rr_{k,0}^{h,sa}(\Omega),\pi)$ is a fiber bundle with fiber $\pi^{-1}(X_0)\simeq\Ff\Rr_{0,0}^{h,sa}(\Omega|_{X_1\times Y_1})$ at $X_0$.    
\end{proof}

\begin{corollary}\label{c:skew-k-path-connected}
	Let $X$ and $Y$ be two Banach space over $\KK$, and $\Omega\colon X\times Y\to\KK$ be a bounded non-degenerate sesquilinear form. Then the space $\Ff\Rr_{k,0}^{h,sa}(\Omega)$ is path connected for each non-negative integer $k$.
\end{corollary}

\begin{proof}
	By Lemma \cite[A.4.4]{BoZh18}, the space $G(k,X)$ is path connected. By Lemma \ref{l:transversal-real-skew}, the space $\Ff\Rr_{0,0}^{h,sa}(\Omega|_{X_1\times Y_1})$ is path connected. By Lemma \ref{l:skew-adjoint-k-fibration}, the space $\Ff\Rr_{k,0}^{h,sa}(\Omega)$ is path connected.
\end{proof}

\begin{proof}[Proof of Theorem \ref{t:path-components-real-skew}]
	Let $a\in\{0,1\}$. By Lemma \ref{l:reduce-to-transversal-case}.a and Proposition \ref{p:real-skew-adjoint-Hilbert}, the set $\Ff\Rr_{a;0}^{-1,sa}(\Omega)$ is nonempty. Take a  $T\in\Ff\Rr_{a;0}^{-1,sa}(\Omega)$ and use the notations in Lemma \ref{l:reduce-to-transversal-case}. By Lemma \ref{l:reduce-to-transversal-case}, we have $T=T_0\oplus T_1$,
	where
	\begin{align*}
		T_0\ :&=X_0\times\{0\}\in\Ff\Rr^{-1,sa}(\Omega|_{X_0\times Y_0}),\\
		T_1\ :&=T\cap(X_1\times Y_1)\in\Ff\Rr_{0,0}^{-1.sa}(\Omega|_{X_1\times Y_1}).	    	
	\end{align*} 
	Then we have 
	\begin{align}\label{e:T0=a mod2}
		a=k=\dim X_0\mod 2.
	\end{align}
	By Proposition \ref{p:real-skew-adjoint-Hilbert}, there is a path $c_0\colon[0,1]\to \Ff\Rr^{-1.sa}(\Omega|_{X_0\times Y_0})$ with $c_0(0)=T_0$ and $ c_0(1)\in\Ff\Rr_{a,0}^{-1,sa}(\Omega|_{X_0\times Y_0})$. Define the path 	
	$c_1\colon[0,1]\to \Ff\Rr^{-1.sa}(\Omega)$ by
	$c_1(t)=c_0(t)\oplus T_1$ for $t\in[0,1]$. Then we have $c_1(0)=T$ and $ c_1(1)\in\Ff\Rr_{a,0}^{-1,sa}(\Omega)$. By Corollary \ref{c:skew-k-path-connected}, the space 
	$\Ff\Rr_{a,0}^{-1,sa}(\Omega)$ is path connected.
	Then the space $\Ff\Rr_{a;0}^{-1,sa}(\Omega)$ is path connected. 
	
	By Theorem \ref{t:stable-mod-2-index}, there is no path joining $\Ff\Rr_{0;0}^{-1.sa}(\Omega)$ and $\Ff\Rr_{1;0}^{-1.sa}(\Omega)$. Then the space $\Ff\Rr^{-1.sa}(\Omega)$ has exactly two path components $\Ff\Rr_{0;0}^{-1.sa}(\Omega)$ and $\Ff\Rr_{1;0}^{-1.sa}(\Omega)$.		    
\end{proof}


\bibliography{Hamiltonian}
\bibliographystyle{plain}

\end{document}